\documentclass[reqno,12pt]{amsart}

\usepackage[T1]{fontenc}
\usepackage[utf8]{inputenc}
\usepackage{lmodern}
\usepackage{amssymb,mathtools}
\usepackage{mathrsfs}
\usepackage{amscd}
\usepackage{graphicx}
\usepackage{array}
\usepackage{booktabs}
\usepackage{multirow}
\usepackage{xcolor}
\usepackage{microtype}
\usepackage[a4paper,margin=1in]{geometry}
\allowdisplaybreaks[4]
\theoremstyle{plain}
\newtheorem{theorem}{Theorem}[section]
\newtheorem{proposition}[theorem]{Proposition}

\newtheorem{lemma}[theorem]{Lemma}

\theoremstyle{definition}

\theoremstyle{remark}
\newtheorem{remark}[theorem]{Remark}

\numberwithin{equation}{section}

\newcommand{\Om}{\Omega}
\newcommand{\Ga}{\Gamma}
\newcommand{\Q}{Q}

\newcommand{\R}{\mathbb{R}}
\newcommand{\dd}{\,\mathrm{d}}

\usepackage[
  colorlinks=true,
  linkcolor=blue,
  citecolor=green,
  urlcolor=blue
]{hyperref}

\title[Inverse Kelvin--Voigt equations]{Simultaneous recovery of elastic and viscous coefficients in a Kelvin--Voigt equation}

\author[Qiu]{Dong Qiu}
\address{School of Mathematical Sciences, Zhejiang University
}
\email{qiudong@zju.edu.cn}

\author[S. Fu]{Song-Ren Fu}
\address{School of Mathematics and Statistics, Northwestern Polytechnical University}
\email{songrenfu@nwpu.edu.cn}

\author[Ye]{Yeqiong Ye}
\address{School of Mathematical Sciences, Zhejiang University}
\email{yeyeqiong@zju.edu.cn}

\author[Zhou]{Ting Zhou}
\address{School of Mathematical Sciences, Zhejiang University
}
\email{ting\_zhou@zju.edu.cn}

\date{}

\begin{document}

\begin{abstract}
In this paper, we study the simultaneous recovery of a spatially dependent elastic coefficient and a time- and space-dependent viscosity coefficient in a scalar Kelvin--Voigt equation. We prove that the displacement-to-traction
map over any fixed time interval uniquely determines both coefficients, provided that the Riemannian metric induced by the viscosity is simple at each time. In the velocity formulation, the elastic term gives rise to a Volterra memory term involving second-order spatial derivatives,
with a corresponding memory contribution to the boundary traction.
We first determine all normal derivatives of both coefficients at the
boundary using highly oscillatory boundary data. We then construct exact
parabolic exponential solutions in the presence of memory and determine
the viscosity by boundary rigidity. With the viscosity determined, an
approximate identity argument reduces recovery of the elastic coefficient
to injectivity of the scalar geodesic ray transform. Finally, numerical
experiments illustrate a regularized finite-dimensional reconstruction
of both coefficients from finitely many moments of the
displacement-to-traction map.
\end{abstract}

\subjclass[2020]{35R30, 35K20}

\keywords{Kelvin--Voigt equation, inverse boundary value problem, Volterra memory, boundary rigidity, geodesic ray transform}

\maketitle

\setlength{\parskip}{0.8ex} 

\section{Introduction}
Viscoelastic solids exhibit both a recoverable deformation response and energy dissipation during motion. The Kelvin--Voigt constitutive law models these effects by adding a viscous stress proportional to the strain rate to an elastic stress proportional to the strain \cite{Christensen1982,SerianiOliveira2020}. This rate-dependent formulation has been used to study material damping in biomimetic scale-covered beams \cite{EbrahimiEtAl2023}. Kelvin--Voigt elements have also been used to represent the viscoelastic membrane and cytoplasm within a hierarchical model of cellular rheology \cite{HangKangXuGao2021}.

Motivated by applications to heterogeneous materials, we allow the elastic coefficient to depend on position and the viscosity to depend on both position and time. In a scalar displacement approximation, the resulting stress flux is $a(x)\nabla u+b(t,x)\nabla\partial_t u$. Combining the divergence of this flux with the balance of linear momentum, and normalizing the mass density to one, yields the equation \eqref{eq:model}. Physically, predicting the response of a heterogeneous viscoelastic body requires knowledge of its internal stiffness and damping. These parameters are often difficult to measure directly, whereas boundary displacements can be prescribed and the resulting tractions observed. This motivates us to study the simultaneous recovery of these coefficients from suitable boundary measurements. 

Let $T>0$ and let $\Om\subset\R^n$, $n\ge2$, be a connected bounded domain with smooth boundary $\partial\Omega=\Ga$.
Throughout this paper, we set
\begin{equation*}
 Q=(0,T)\times\Om,\quad \Sigma=(0,T)\times\Ga.
\end{equation*}
We consider the following scalar Kelvin--Voigt equation
\begin{equation}\label{eq:model}
\begin{cases}
\partial_t^2 u-\operatorname{div}(a(x)\nabla u)-\operatorname{div}(b(t,x)\nabla \partial_tu)=0&\text{in }Q,\\
u=f&\text{on }\Sigma,\\
u(0,\cdot)=\partial_tu(0,\cdot)=0&\text{in }\Om.
\end{cases}
\end{equation}
Here $u$ denotes the scalar displacement, $a(x)$ is the elastic coefficient and $b(t,x)$ is the time-dependent viscosity coefficient. The terms $a\nabla u$ and $b\nabla\partial_tu$
are the elastic and viscous stress fluxes, respectively. We note that the mass density of a given medium has been normalized to one.

Henceforth, we assume that the coefficients $a\in C^\infty(\overline\Omega)$ and $b\in C^\infty([0,T]\times\overline\Omega)$ are real-valued and that there exist positive constants $a_*,b_*$ such that
\begin{equation}\label{eq:positivity}
a(x)\ge a_*>0,\quad b(t,x)\ge b_*>0.
\end{equation}

\subsection{Formulation of the inverse problem and main result}
We next formulate the considered inverse problem. 
Define the stress flux by
\begin{equation*}
\sigma(u)\vcentcolon=a\nabla u+b\nabla \partial_tu.
\end{equation*}
For the moment, we formally introduce the displacement-to-traction map by
\begin{equation}
\Lambda_{a,b}^T f\vcentcolon=\langle\sigma(u),\nu\rangle|_\Sigma=(a\partial_\nu u+b\partial_\nu \partial_tu)|_\Sigma.
\label{eq:response}
\end{equation}
Here, $\nu$ denotes the outward unit normal vector field of $\Omega$ along $\Ga$. The well-posedness of the forward problem \eqref{eq:model} is established in Proposition~\ref{pro:well-posedness-f=K+h}, so that the map $\Lambda_{a,b}^T f$ is well-defined. We note that, for real smooth solutions, we have the energy identity
\begin{equation*}
\frac{\dd}{\dd t}\frac12\int_\Om (|\partial_tu|^2+a|\nabla u|^2)\dd x+\int_\Om b|\nabla \partial_tu|^2\dd x
 =\int_\Ga \sigma(u)\cdot\nu\,\partial_tu\dd S.
\end{equation*}
This provides a natural derivation of the displacement-to-traction map.

Let $K_+$ and $K_-$ be the forward and backward Volterra operators, 
\begin{equation*}
 (K_+h)(t)=\int_0^t h(s)\dd s,\quad
 (K_-h)(t)=\int_t^T h(s)\dd s.
\end{equation*}
We consider lateral Dirichlet data in the class
\begin{equation*}
\mathcal D_T:=K_+C_c^\infty(\Sigma)=\{K_+h:h\in C_c^\infty(\Sigma)\}.
\end{equation*}
Thus, for $f\in\mathcal D_T$, the boundary velocity
$h=\partial_t f$ is smooth and compactly supported in $\Sigma$.
The displacement $f$ vanishes near $t=0$ but is constant, and need
not vanish, near $t=T$. Section~\ref{sec:forward} establishes
well-posedness for a larger Sobolev class of boundary data and the
boundedness of the associated displacement-to-traction map.

The concerned inverse problem is stated as follows.

{\bf (IP)} \emph{Can one simultaneously recover the elastic coefficient $a(x)$ and the time-dependent viscosity $b(t,x)$ by the boundary displacement-to-traction map $\Lambda_{a,b}^T$?}


For fixed time $t$, the relevant Riemannian metric associated with $b$ is given by
\begin{equation}\label{eq:metric}
g_b(t)=b^{-1}(t,\cdot)e,\quad e=\dd s^2=\dd x_1^2+\cdots \dd x_n^2,
\end{equation}
where $e$ denotes the standard Euclidean metric.
A compact connected Riemannian manifold $(\mathcal M,g)$ with smooth boundary is \emph{simple} if its boundary is strictly convex, it is nontrapping, and it has no conjugate points.
In a simple manifold, any two points are joined by a unique geodesic, which is minimizing (see~\cite[Section~3.8]{PaternainSaloUhlmann2023}).

Below, we give an affirmative answer to the considered inverse problem.
\begin{theorem}\label{thm:main}
Let $T>0$ be given. For $j=1,2$, let $a_j\in C^\infty(\overline\Om)$ and $b_j\in C^\infty([0,T]\times\overline\Om)$ be real-valued and satisfy the condition \eqref{eq:positivity}.
Assume that the geometry $(\overline\Om,g_{b_j}(t))$ is simple for every $t\in[0,T]$.
If
\begin{equation}\label{eq:dataequality}
\Lambda_{a_1,b_1}^T f=\Lambda_{a_2,b_2}^T f\quad\text{for all }f\in\mathcal D_T,
\end{equation}
then $a_1=a_2$ in $\overline\Omega$ and $b_1=b_2$ in $[0,T]\times\overline\Omega$.
\end{theorem}

\begin{remark}
Theorem~\ref{thm:main} uses the displacement-to-traction map on the whole lateral boundary for all $f\in\mathcal D_T$, over any fixed $T>0$. 
Simplicity is required only for $g_{b_j}(t)=b_j^{-1}(t,\cdot)e$, with no geometric condition on $a_j^{-1}e$.
\end{remark}

\begin{remark}
The positivity of $a$ in \eqref{eq:positivity} is assumed only for the physical interpretation of the energy. Mathematically, Theorem \ref{thm:main} holds for any real-valued $a \in C^{\infty}(\overline{\Omega})$.
\end{remark}

\subsection{Related literature and main difficulties}

Inverse problems for dissipative evolution equations have attracted considerable attention over the past decades, motivated by applications such as parameter identification, boundary data completion, and source or coefficient determination, as well as by the mathematical challenges arising from the loss of time reversibility and the smoothing effect of the associated semigroup. In particular, strongly damped wave and Kelvin--Voigt equations have been studied with a variety of unknowns and measurement schemes. Carleman estimates have been used to obtain inverse source results for Kelvin--Voigt systems~\cite{ImanuvilovYamamoto2019} and to recover a lower-order coefficient in a strongly damped wave equation~\cite{Wu2012}.

Li and Zhang~\cite{LiZhang2024} determined, up to isometry, the metric shared by the elastic and damping terms of a strongly damped wave equation. Related work in nonlinear acoustics addresses the recovery of a time-dependent potential and a nonlinear coefficient \cite{LiZhangAcoustics2024}, as well as the simultaneous recovery of wave speed, diffusivity, and nonlinearity from boundary Cauchy data at the first and second harmonics under time-harmonic boundary excitation at a sufficiently high frequency~\cite{AcostaPalacios2026}. Yamamoto et al. \cite{YamamotoZhangZhao2026} established H\"older stability in recovering five spatially varying coefficients in a three-dimensional viscoelastic system, using finitely many measurements along with fixed-time interior data. Jiang et al.~\cite{JiangLinNakamura2026} proved conditional uniqueness for the shear wave speed in two isotropic viscoelastic models from a single interior displacement record. 

For equations with memory, boundary measurements have been used to
recover memory kernels, including from partial boundary
data~\cite{BukhgeimDyatlovIsakov1997,BukhgeimDyatlovUhlmann2002}, and
lower-order coefficients in wave equations with nonlocal
attenuation~\cite{KianRezigUhlmann2026}. Memory kernels in strongly
damped wave equations have also been recovered from supplementary
observations, including a flux
condition~\cite{Colombo2007,ColomboGuidetti2009}.
In the parabolic setting, coefficient recovery from boundary
measurements was studied in~\cite{CanutoKavian2001}. Carleman estimates
have also yielded stability results for equations with second-order
spatial memory~\cite{BaraniBalanSakthivelBalachandran2008}.

For the viscoacoustic wave equation, exterior source-to-solution data
have been shown, under suitable geometric assumptions, to determine
the sound speed and all time derivatives at zero of the kernel of a
memory term involving second-order spatial
derivatives~\cite{CoviDeHoopSalo2026}.
High-frequency propagation in the presence of weakly singular hereditary
memory has been studied in~\cite{GrebnevStefanov2026}.

Our approach to recovering the time-dependent viscosity draws on the geometric approach to inverse boundary problems for nonautonomous heat equations developed in Feizmohammadi's work~\cite{Feizmohammadi2024}. In the velocity
formulation of the Kelvin--Voigt equation, the elastic term gives rise
to a Volterra memory term involving two spatial derivatives. The boundary traction contains a corresponding memory contribution that must also be controlled in the inverse analysis.

We first construct boundary-localized, highly oscillatory solutions and use them to determine all normal derivatives of both coefficients at the boundary. We then construct exact parabolic exponential solutions by showing that the conjugated memory operator, composed with the inverse of the conjugated parabolic operator, is small as the
exponential parameter tends to infinity. After controlling the memory contribution to the boundary traction, we
recover the viscosity by boundary rigidity. Once the viscosity is known, a strong approximate-identity argument reduces the recovery of the elastic coefficient to the injectivity of the scalar geodesic ray transform applied to $(a_1-a_2)/b(t,\cdot)$.

The remainder of the paper is organized as follows. In Section~\ref{sec:forward}, we develop the parabolic formulation of the original equation, establish well-posedness using a global exponential weight, and derive two useful identities related to the recovery of coefficients. In Section~\ref{sec:boundary}, we study the boundary determination of all normal derivatives of both coefficients at the boundary. Section~\ref{sec:exponential} focuses on constructing exact exponential solutions with geodesic distance functions as phases in the presence of the Volterra memory term. The proof of the main theorem will be provided in Section~\ref{sec:proof-theorem}. Finally, Section~\ref{sec:numerics} presents a finite-dimensional reconstruction method and numerical experiments.

\section{The forward problem and boundary identities}\label{sec:forward}
In this section, we begin by reformulating the Kelvin--Voigt equation as a parabolic
equation with a Volterra memory term, and then establish well-posedness
and the boundary identities used in the inverse analysis. Let $v=\partial_t u$. Since $u(0,\cdot)=0$,
we have $u=K_+v$, where $K_+$ acts in the time variable. Set
\begin{equation*}
D_a=\operatorname{div}(a\nabla),\quad
P_b^\pm=\pm\partial_t-\operatorname{div}(b\nabla),\quad
L_{a,b}^\pm=P_b^\pm-K_\pm D_a.
\end{equation*}
As $a$ is independent of time, $D_aK_+v=K_+D_av$.
Hence \eqref{eq:model} is equivalent to
\begin{equation}
\begin{cases}
L_{a,b}^+v=0 & \text{in }Q,\\
v=h=\partial_t f & \text{on }\Sigma,\\
v(0,\cdot)=0 & \text{in }\Omega,
\end{cases}
\qquad u=K_+v.
\label{eq:velocity}
\end{equation}

\subsection{Trace spaces and well-posedness results}

We introduce the following Hilbert spaces
\begin{align*}
\mathbb Y&=L^2(Q),&
\mathbb W&=H^1(0,T;L^2(\Om))
          \cap L^2(0,T;H^2(\Om)),\\
\mathbb X&=\mathbb W\cap L^2(0,T;H_0^1(\Om)),&
\mathbb X_+&=\{z\in\mathbb X:z(0)=0\},\quad
\mathbb X_-=\{z\in\mathbb X:z(T)=0\}.
\end{align*}
The spaces $\mathbb X$ and $\mathbb X_\pm$ are equipped with the norm
inherited from $\mathbb W$. The embedding
$\mathbb W\hookrightarrow C([0,T];H^1(\Om))$ ensures that the endpoint
conditions are well defined. The lateral Dirichlet and normal trace spaces are
\begin{align}
\mathbb H&=H^{3/4}(0,T;L^2(\Ga))\cap L^2(0,T;H^{3/2}(\Ga)),\notag\\
\mathbb H_+&=\{h\in\mathbb H:h(0)=0\},\quad
\mathbb H_-=\{k\in\mathbb H:k(T)=0\},\notag\\
\mathbb N&=H^{1/4}(0,T;L^2(\Ga))\cap L^2(0,T;H^{1/2}(\Ga)).
\label{eq:spaces}
\end{align}
The endpoint conditions defining $\mathbb H_\pm$ are understood in
$H^{1/2}(\Ga)$. The joint initial--lateral trace theorem provides
bounded liftings $E_\pm:\mathbb H_\pm\to\mathbb W$ with the prescribed
lateral trace and zero initial or terminal trace, respectively.
Moreover, for every $z\in\mathbb W$, we have
\begin{equation}
\|\partial_\nu z\|_{\mathbb N}\le C\|z\|_{\mathbb W}.
\label{eq:normaltrace}
\end{equation}
We will use the scalar Dirichlet maximal-regularity and trace results
of~\cite[Theorem~2.1 and Section~3]{DenkHieberPruss2007}. The positivity and smoothness of $b$ allow us to apply Dirichlet
maximal regularity. For $\eta\in\mathbb H_\pm$, the lifting $E_\pm\eta\in\mathbb W$ also satisfies
$L_{a,b}^\pm E_\pm\eta\in\mathbb Y$.

\begin{proposition}\label{wellpose-Lab}\label{prop:wp}
For any \(F \in \mathbb{Y}\) and \(h \in \mathbb{H}_+\), the initial-boundary value problem
\begin{equation}\label{eq:forward-mem}
    \begin{cases}
L_{a,b}^+ v = F & \text{in } Q, \\
v = h, & \text{on } \Sigma\\
v(0) = 0 & \text{in } \Omega,
\end{cases}
\end{equation}
admits a unique solution \(v \in \mathbb{W}\), satisfying the estimate
\begin{equation}\label{eq:est-v}
\|v\|_{\mathbb{W}} + \|\partial_\nu v\|_{\mathbb{N}} \leq C\bigl(\|F\|_{\mathbb{Y}} + \|h\|_{\mathbb{H}}\bigr). 
\end{equation}
The terminal problem
\begin{equation*}
L_{a,b}^-p=F\ \text{in }Q,\quad
p=k\ \text{on }\Sigma,\quad
p(T,\cdot)=0\ \text{in }\Omega
\end{equation*}
is also well posed for $F\in\mathbb Y$ and $k\in\mathbb H_-$. In particular, the map $L_{a,b}^-:\mathbb X_-\to\mathbb Y$ is an isomorphism.
\end{proposition}

\begin{proof}
We begin by deriving a uniform \emph{a priori} estimate for the shifted heat equation \((P_b^+ + \lambda_0) z = g\) for any fixed \(\lambda_0 \geq 1\).
By \cite[Theorem~2.1]{DenkHieberPruss2007}, this equation admits a unique solution $z \in \mathbb{X}_+$ for any $g \in \mathbb Y$.
We first derive the energy identities for smooth data. The general case $g\in\mathbb Y$ follows by approximation and maximal regularity. Testing $(P_b^++\lambda_0)z=g$ with $\bar z$, integrating over $Q$, and taking real parts yield
\begin{equation*}
\frac12\|z(T)\|_{L^2(\Omega)}^2+\int_Q b|\nabla z|^2\,\dd x\dd t+\lambda_0\|z\|_{\mathbb Y}^2=\operatorname{Re}\int_Q g\bar z\,\dd x\dd t.
\end{equation*}
Dropping the nonnegative endpoint term and using $b\ge b_*>0$
and the Cauchy--Schwarz inequality, we obtain
\begin{equation*}
\|z\|_{\mathbb Y}\le \lambda_0^{-1}\|g\|_{\mathbb Y},\quad \|\nabla z\|_{\mathbb Y}\le b_*^{-1/2}\lambda_0^{-1/2}\|g\|_{\mathbb Y}.
\end{equation*}

We next test by \(\bar{z}_t\) to obtain
\begin{equation*}
\|z_t\|_{\mathbb{Y}}^2 + \frac{1}{2}\int_\Omega \big(b|\nabla z|^2 + \lambda_0|z|^2 \,dx\big)\Bigr|_0^T = \operatorname{Re} \int_Q g\bar{z}_t \,\dd x \dd t + \frac{1}{2}\int_Q b_t|\nabla z|^2 \,\dd x \dd t.
\end{equation*}
Using $z(0)=0$ in $H^1(\Om)$, we can get
\begin{equation*}
\|z_t\|_{\mathbb Y}^2\le \|g\|_{\mathbb Y}\|z_t\|_{\mathbb Y} +\frac12\|b_t\|_{L^\infty(Q)}\|\nabla z\|_{\mathbb Y}^2.
\end{equation*}
By Young's inequality and the estimate for $\nabla z$, we have $\|z_t\|_{\mathbb Y}\le C\|g\|_{\mathbb Y}$ uniformly for $\lambda_0\ge1$. Regarding \(-\operatorname{div}(b\nabla z) = g - z_t - \lambda_0 z\) as a stationary elliptic problem for a.e.\ \(t\), and invoking Dirichlet elliptic regularity, we conclude that
\begin{align}\label{eq:est-z}
\|z\|_{\mathbb{W}} + \lambda_0\|z\|_{\mathbb{Y}} + \sqrt{\lambda_0}\|\nabla z\|_{\mathbb{Y}} \leq C\|g\|_{\mathbb{Y}}, \quad \lambda_0 \geq 1,
\end{align}
where the constant \(C > 0\) is independent of \(\lambda_0\).

We now turn to the inhomogeneous Volterra problem. Using the trace lifting operator, we set \(v = \tilde{h} + w\) with $\tilde{h}=E_+ h$ and \(w \in \mathbb{X}_+\), which reduces the equation \eqref{eq:forward-mem} to
\begin{equation}\label{eq:lifting}
    P_b^+ w - K_+ D_a w = \tilde{F},
\end{equation}
where \(\tilde{F} = F - P_b^+ \tilde{h} + K_+ D_a \tilde{h}\) and satisfies
\begin{equation}\label{eq:est-tildeF}
    \|\tilde{F}\|_{\mathbb{Y}} \leq C\bigl(\|F\|_{\mathbb{Y}} + \|h\|_{\mathbb{H}}\bigr).
\end{equation}
Let \(w = e^{\lambda_0 t} z\). We see that
\begin{equation}\label{eq:z-V}
    (P_b^+ + \lambda_0) z - V_{\lambda_0} D_a z = e^{-\lambda_0 t}\tilde{F},
\end{equation}
where \((V_{\lambda_0} q)(t) = \int_0^t e^{-\lambda_0(t-s)} q(s) \,ds\).
Let $S_{\lambda_0}=(P_b^++\lambda_0)^{-1}$ with homogeneous initial and lateral data. The relevant bounded maps are
\begin{equation*}
D_a:\mathbb X_+\longrightarrow\mathbb Y,\quad V_{\lambda_0}:\mathbb Y\longrightarrow\mathbb Y,\quad
S_{\lambda_0}:\mathbb Y\longrightarrow\mathbb X_+.
\end{equation*}
Estimate~\eqref{eq:est-z} shows that $S_{\lambda_0}:\mathbb Y\to\mathbb X_+$ is bounded uniformly for $\lambda_0\ge1$.

Given $z \in \mathbb{X}_+$, we consider the following problem:
\begin{equation}\label{eq:tilde-z}
    \begin{cases}
(P_b^+ + \lambda_0) \tilde{z} = V_{\lambda_0} D_a z+ e^{-\lambda_0 t}\tilde{F} & \text{in } Q,\\
\tilde{z}= 0, & \text{on } \Sigma\\
\tilde{z}(0) = 0 & \text{in } \Omega.
\end{cases}
\end{equation}
Since $\|V_{\lambda_0}\|_{\mathbb Y\to\mathbb Y} \le\lambda_0^{-1}$ and $D_a:H^2(\Om)\to L^2(\Om)$ is bounded, we have $V_{\lambda_0} D_a z+ e^{-\lambda_0 t}\tilde{F} \in \mathbb{Y}$.
By maximal regularity for $P_b^++\lambda_0$, problem~\eqref{eq:tilde-z} has a unique solution $\tilde z\in\mathbb X_+$.
Moreover, estimate~\eqref{eq:est-z} gives
\begin{align}\label{eq:est-tildez}
\|\tilde{z}\|_{\mathbb X_+} \leq C \|V_{\lambda_0} D_a z+ e^{-\lambda_0 t}\tilde{F}\|_{\mathbb Y} \leq C C_a\lambda_0^{-1}\|z\|_{\mathbb{X}_+} + C\|\tilde{F}\|_{\mathbb{Y}}.
\end{align}
Define $\mathfrak T:\mathbb X_+\to\mathbb X_+$ by
$\mathfrak Tz=\tilde z$, where $\tilde z$ is the solution
of~\eqref{eq:tilde-z}. For $z_1,z_2\in\mathbb X_+$,
\begin{equation*}
\|\mathfrak Tz_1-\mathfrak Tz_2\|_{\mathbb X_+}\le C\|V_{\lambda_0} D_a(z_1-z_2)\|_{\mathbb Y}\le CC_a\lambda_0^{-1}\|z_1-z_2\|_{\mathbb X_+},
\end{equation*}
where $C_a$ is the operator norm of
$D_a:\mathbb X_+\to\mathbb Y$.
Choose $\lambda_0\ge1$ so that $CC_a/\lambda_0<1$.
Banach's fixed-point theorem then gives a unique
$z\in\mathbb X_+$ solving~\eqref{eq:z-V}.
At the fixed point, estimate~\eqref{eq:est-tildez} allows the term
involving $\|z\|_{\mathbb X_+}$ to be absorbed.
Returning to $v=\tilde h+e^{\lambda_0 t}z$ and
using~\eqref{eq:est-tildeF} yields~\eqref{eq:est-v}.
The terminal result follows by time reversal.
\end{proof}

\begin{proposition}\label{pro:well-posedness-f=K+h}
Let $h\in\mathbb H_+$ and set $f=K_+h$. Then $v\in\mathbb W$
solves~\eqref{eq:forward-mem} with $F=0$ if and only if
$u=K_+v$ solves~\eqref{eq:model} and belongs to $H^2(0,T;L^2(\Omega))\cap H^1(0,T;H^2(\Omega)).$
Thus the displacement-to-traction map induces a bounded
boundary-velocity-to-traction map
\begin{equation}
\mathcal R_{a,b}^+:\mathbb H_+\rightarrow\mathbb N,\quad h\longmapsto\Lambda_{a,b}^T(K_+h)
=b\,\partial_\nu v+aK_+\partial_\nu v.
\label{eq:velocityresponse}
\end{equation}
\end{proposition}

\begin{proof}
If $v\in\mathbb W$, then $u=K_+v\in H^2(0,T;L^2(\Omega))\cap H^1(0,T;H^2(\Omega))$.
Since $a$ is independent of time, $D_au=K_+D_av$. Together with the initial and lateral traces, these identities show
that a solution $v$ of~\eqref{eq:forward-mem} with $F=0$ yields a solution $u$ of~\eqref{eq:model}.

Hence, if $u$ belongs to the stated space and solves~\eqref{eq:model}, then $v=u_t\in\mathbb W$ satisfies
\eqref{eq:forward-mem} with $F=0$ and boundary data $h$. Finally, $\partial_\nu K_+v=K_+\partial_\nu v$. As $K_+$ and multiplication by the smooth boundary coefficients are bounded on $\mathbb N$, the estimate~\eqref{eq:est-v} gives
\begin{equation*}
\|\mathcal R_{a,b}^+h\|_{\mathbb N}\le C\|h\|_{\mathbb H},
\end{equation*}
which proves that $\mathcal R_{a,b}^+:\mathbb H_+\to\mathbb N$ is bounded.
\end{proof}

\subsection{Boundary and volume integral identities}
Define the boundary operators
\begin{equation*}
\mathcal C_{a,b}^\pm z
=b\,\partial_\nu z+aK_\pm\partial_\nu z.
\end{equation*}
For $k\in\mathbb H_-$, let $p$ solve $L_{a,b}^-p=0$ in $Q$,
with $p=k$ on $\Sigma$ and $p(T,\cdot)=0$. We define the auxiliary
terminal boundary map by
$\mathcal R_{a,b}^-k=\mathcal C_{a,b}^-p$.
The pairings in Green's formula and the inverse identities are complex bilinear. Green's formula gives, for $v,p\in\mathbb W$,
\begin{equation}
\int_Q\big((L_{a,b}^+v)p-vL_{a,b}^-p\big)\,\dd x\dd t=\int_\Om vp\,\dd x\Big|_0^T+\int_\Sigma\big(v\mathcal C_{a,b}^-p
-(\mathcal C_{a,b}^+v)p\big)\,\dd S\dd t.
\label{eq:green}
\end{equation}
For the Volterra terms, we use $K_+^*=K_-$, Fubini's theorem,
and the spatial Green formula for $D_a$. The time independence of
$a$ ensures that $D_a$ commutes with $K_\pm$.
The identity follows first for smooth functions and extends to
$\mathbb W$ by approximation and the trace estimates. We have the following two useful integral identities.

\begin{lemma}\label{lem:identities}
Let $h,k\in C_c^\infty(\Sigma)$, and let $v_1,p_2\in\mathbb W$ be the corresponding solutions of
\[
L_{a_1,b_1}^+v_1=0,\quad v_1|_\Sigma=h,\quad v_1(0,\cdot)=0,\quad
L_{a_2,b_2}^-p_2=0,\quad p_2|_\Sigma=k,\quad p_2(T,\cdot)=0.
\]
Then, under the condition~\eqref{eq:dataequality},
\begin{equation}
\mathcal B(v_1,p_2):=\int_\Sigma\bigl((\mathcal C_{a_1,b_1}^+v_1)p_2-v_1\mathcal C_{a_2,b_2}^-p_2\bigr)\,\dd S\dd t=0
\label{eq:boundaryidentity}
\end{equation}
and
\begin{equation}
\int_Q(b_1-b_2)\nabla v_1\cdot\nabla p_2\,\dd x\dd t+\int_Q(a_1-a_2)\nabla(K_+v_1)\cdot\nabla p_2\dd x\dd t=0.
\label{eq:mixedidentity}
\end{equation}
\end{lemma}

\begin{proof}
For either coefficient pair, applying~\eqref{eq:green} to the
corresponding forward and terminal solutions gives
\begin{equation*}
\int_\Sigma(\mathcal R_{a_j,b_j}^+h)k\,\dd S\dd t=\int_\Sigma h(\mathcal R_{a_j,b_j}^-k)\,\dd S\dd t, \quad j=1,2,
\end{equation*}
because the endpoint term vanishes. Since $K_+h\in\mathcal D_T$, assumption~\eqref{eq:dataequality} implies
$\mathcal R_{a_1,b_1}^+h=\mathcal R_{a_2,b_2}^+h$. The transposition identity for the second coefficient pair therefore
yields
\begin{equation*}
\mathcal B(v_1,p_2)=\int_\Sigma \bigl((\mathcal R_{a_1,b_1}^+-\mathcal R_{a_2,b_2}^+)h\bigr)k\,\dd S\dd t=0.
\end{equation*}

Next, we apply~\eqref{eq:green} for $(a_1,b_1)$ to $v_1$ and $p_2$, and use $L_{a_2,b_2}^-p_2=0$. Integration by parts and
$K_+^*=K_-$ give
\begin{equation*}
\mathcal B(v_1,p_2)=\int_Q(b_1-b_2)\nabla v_1\cdot\nabla p_2\,\dd x\dd t+\int_Q(a_1-a_2)\nabla(K_+v_1)\cdot\nabla p_2\,\dd x\dd t.
\end{equation*}
Here the time independence of $a_1-a_2$ allows it to pass through $K_\pm$. Combining this equality with
$\mathcal B(v_1,p_2)=0$ proves~\eqref{eq:mixedidentity}, as desired. 
\end{proof}

\begin{remark}
We will use the boundary identity~\eqref{eq:boundaryidentity} to recover $b$ without prior knowledge of $a$ in the interior. Once $b$ is known, the volume identity~\eqref{eq:mixedidentity} yields the recovery of $a$. The proof uses only measurements with boundary inputs $f=K_+h$, $h\in C_c^\infty(\Sigma)$, and hence requires $\Lambda_{a,b}^T$ only on the class $\mathcal D_T$
specified in Theorem~\ref{thm:main}.
\end{remark}

\section{Boundary determination of the coefficients}\label{sec:boundary}
In this section, we determine all normal derivatives of the coefficients $a$ and $b$ on the boundary from the displacement-to-traction map through the construction of boundary-localized, highly oscillatory solutions. Throughout this section, all equations and measurement maps are complexified, and pairings of complex spatial covectors are bilinear.

We now state the following boundary determination result for $a$ and $b$.
\begin{proposition}\label{prop:boundary}
Let $(a_j,b_j)$, $j=1,2$, satisfy the conditions of Theorem~\ref{thm:main}, and suppose that \eqref{eq:dataequality} holds. Then, for every $m\geq0$,
\begin{equation}
\partial_\nu^m a_1=\partial_\nu^m a_2\quad\text{on }\Ga,\quad \partial_\nu^m b_1=\partial_\nu^m b_2\quad\text{on }[0,T]\times\Ga.
 \label{eq:joint-boundary-jets}
\end{equation}
\end{proposition}

\subsection{Construction of the phase and amplitude functions}

Fix $(t_0,x_0)\in(0,T)\times\Gamma$, and let $U$ be a sufficiently small neighborhood of $x_0$. We choose Euclidean boundary normal coordinates $y=(y_1,y')$ on $U$ such that $y_1=0$ on $\Gamma\cap U$ and $y_1<0$ in $\Omega\cap U$.
Let $U'$ denote the $y'$-coordinate image of $\Gamma\cap U$. In these coordinates,
\[
\nu=\partial_{y_1}\quad\text{on }\Gamma\cap U,\quad e=\dd y_1^2+e_{\alpha\beta}(y)\,\dd y_\alpha\dd y_\beta.
\]
Thus $\rho=-y_1>0$ is the inward normal distance in $\Omega\cap U$.

We next construct an approximate geometric optics solution localized near the boundary, taking the following WKB form
\begin{equation}
 v_\tau^{\mathrm{app}}=e^{i\tau t+\tau\Phi}\vartheta w_\tau^{\rm amp},
 \quad w_\tau^{\rm amp}=\sum_{k=0}^{N-1}\tau^{-k}w_k,
 \label{eq:bd-ansatz}
\end{equation}
where \(\tau > 0\) is a large parameter, \(\Phi(y)\) is a time-independent complex phase function, and $\vartheta\in C_c^\infty(\mathbb R)$ is a cutoff function in the normal variable $y_1$, equal to $1$ for $|y_1|\le\varepsilon/4$
and vanishing for $|y_1|\ge\varepsilon/2$. The amplitudes \(w_k\) satisfy the same boundary condition $w_k|_{\Gamma} = \delta_{k0}\chi,$ where \(\chi \in C_c^\infty((0,T) \times U')\). 

Direct calculation of the conjugated operator $\mathcal{P}_\tau w = e^{-i\tau t - \tau \Phi} L_{a,b}^+(e^{i\tau t + \tau \Phi}w)$ yields
\begin{align}
 \mathcal{P}_\tau w = -\tau^2 b E_\Phi w - \tau T_b w + P_b^+ w - \mathcal{K}_\tau(\tau^2 a E_\Phi w + \tau S_a w + D_a w),
 \label{eq:bd-conjugation}
\end{align}
where the eikonal defect is \(E_\Phi := \langle d\Phi, d\Phi\rangle_e\), and the transport and memory operators are defined as
\begin{align}
 T_bw&=2b\langle d\Phi,dw\rangle_e
       +(\langle db,d\Phi\rangle_e+b\Delta\Phi-i)w,
       \label{eq:bd-Tb}\\
 S_aw&=2a\langle d\Phi,dw\rangle_e
       +(\langle da,d\Phi\rangle_e+a\Delta\Phi)w,
       \label{eq:bd-Daphi}\\
 \mathcal K_\tau F&=e^{-i\tau t}K_+(e^{i\tau t}F).
 \label{eq:bd-history-definition}
\end{align}

To eliminate the terms of order \(\tau^2\) in \eqref{eq:bd-conjugation}, we construct the phase function \(\Phi\) of the following form
\[
\Phi(y) = \sum_{j=0}^{N} \frac{y_1^j}{j!} \Phi_j(y').
\]
We require \(E_\Phi\) to vanish up to order \(N - 1\) at the boundary \(y_1 = 0\). For \(j = 0\), we have 
\[
\Phi_1^2 + e^{\alpha\beta}(0, y')\partial_\alpha \Phi_0 \partial_\beta \Phi_0 = 0.
\]
Choose a real-valued smooth function $\psi$ on $U'$ satisfying
$|d\psi|_{e|_\Gamma}^2=1$, and set $\Phi_0=i\psi$ and $\Phi_1=1$. Then $E_\Phi|_\Gamma=0$. For $1\le j\le N-1$, the condition $\partial_{y_1}^jE_\Phi|_\Gamma=0$ uniquely determines $\Phi_{j+1}$ from $\Phi_0,\cdots,\Phi_j$. Taylor's theorem
therefore gives $E_\Phi=O(\rho^N)$. Moreover, since $\operatorname{Re}\Phi=-\rho+O(\rho^2)$, we have
$\operatorname{Re}\Phi\le -c\rho$ for $0\le\rho<\varepsilon$, provided $\varepsilon>0$ is sufficiently small.

Let $X$ be a Hilbert space, and let $F$ be a smooth $X$-valued function that vanishes near $t=0$. Repeated integration by parts
in time yields the following expansion for $\mathcal K_\tau F$:
\begin{equation}
\mathcal K_\tau F=\sum_{r=0}^{N_{\rm time}-1}\tau^{-r-1}c_r\partial_t^rF+\mathcal E_{N_{\rm time},\tau}F,
 \label{eq:bd-history-expansion}
\end{equation}
where $c_r=(-1)^ri^{-r-1}$ and
\[
\mathcal E_{N_{\rm time},\tau}F=(-1)^{N_{\rm time}}(i\tau)^{-N_{\rm time}}\mathcal K_\tau(\partial_t^{N_{\rm time}}F).
\]
A further integration by parts, together with the boundedness of $K_+$ on $L^2(0,T;X)$, yields the remainder estimate
\begin{equation}
 \|\mathcal E_{N_{\rm time},\tau}F\|_{L^2(0,T;X)}
 \le C_T\tau^{-N_{\rm time}-1}
 \bigl(\|\partial_t^{N_{\rm time}}F\|_{L^2(0,T;X)}
       +\|\partial_t^{N_{\rm time}+1}F\|_{L^2(0,T;X)}\bigr).
 \label{eq:bd-history-remainder}
\end{equation}
Substituting the expansion into the memory terms of \eqref{eq:bd-conjugation}, we obtain the contribution at order \(\tau^{1-k}\) for \(k\ge1\)
\begin{equation}
M_k=-\sum_{r+q=k-1}c_r\partial_t^rS_aw_q-\sum_{r+q=k-2}c_r\partial_t^rD_aw_q.
\label{eq:bd-memory-recursion}
\end{equation}
Since $E_\Phi=O(\rho^N)$, the terms containing $E_\Phi$ in~\eqref{eq:bd-conjugation} have vanishing normal derivatives
on $\Gamma$ through order $N-1$. Matching the remaining powers of $\tau$ gives the following relations for the normal derivatives of the amplitudes at the boundary:
\begin{align}
\partial_\nu^j(T_bw_0)|_\Ga&=0,
&&0\le j\le N-1,\label{eq:prin-amplitude}\\
\partial_\nu^j(-T_bw_k+P_b^+w_{k-1}+M_k)|_\Ga&=0,
&&1\le k\le N-1,\quad0\le j\le N-k-1.
\label{eq:bd-jet-equations}
\end{align}

For $m\ge0$, let
$J_\nu^m c=(\partial_\nu^j c|_\Sigma)_{j=0}^m$
denote the collection of normal derivatives of $c$ up to order $m$, and let $J_\nu^{-1}c$ denote the empty tuple.
Set $\beta=b|_\Sigma$ and $\alpha=a|_\Gamma$. 
For $0\le k\le N-1$ and $0\le\ell\le N-k$, set $W_{k,\ell}=(\partial_\nu^\ell w_k)|_\Sigma.$ 

The following lemma describes how the boundary equations determine the normal derivatives of the amplitudes in terms of those of $a$ and $b$.

\begin{lemma}\label{lem:boundary-amplitude-derivatives}
For $1\le\ell\le N$,
\begin{equation}
W_{0,\ell}=-\frac{\chi}{2\beta}\,\partial_\nu^\ell b+F_{0,\ell}(J_\nu^{\ell-1}b).
\label{eq:W0l}
\end{equation}
For $k,\ell\ge1$ with $k+\ell\le N$,
\begin{equation}
\begin{split}
W_{k,\ell}=\left(-\frac12\right)^{k+1}\frac{\chi}{\beta}\,\partial_\nu^{k+\ell}b+\left(-\frac12\right)^k
\frac{\chi}{i\beta}\,\partial_\nu^{k+\ell-1}a
+F_{k,\ell}(J_\nu^{k+\ell-1}b,J_\nu^{k+\ell-2}a).
\end{split}
\label{eq:Wkl}
\end{equation}
Here, $F_{k,\ell}$ are smooth differential expressions depending strictly on the indicated lower-order normal derivatives of $a,b$, their tangential and time derivatives, and the known functions $\chi,e,\Phi$ and $\beta$.
\end{lemma}

\begin{proof}
We argue by induction on $s=k+\ell$, determining
\[
W_{0,s},\ W_{1,s-1},\ \ldots,\ W_{s-1,1}
\]
in this order. The values $W_{k,0}=\delta_{k0}\chi$ are prescribed. For $k=0$, the equation with $j=\ell-1$
in~\eqref{eq:prin-amplitude} has the leading terms $2\beta W_{0,\ell}+(\partial_\nu^\ell b)\chi.$ All other terms involve lower-order normal derivatives of $b$.
Since $\beta>0$, this determines $W_{0,\ell}$ and
gives~\eqref{eq:W0l}.

For $k\ge1$, the equation with $j=\ell-1$
in~\eqref{eq:bd-jet-equations} contains the terms
\begin{equation}
-2\beta W_{k,\ell}-\beta W_{k-1,\ell+1}.
\label{eq:induc-k}
\end{equation}
The second amplitude derivative is already known at the current induction step. In the two sums defining $M_k$, the indices satisfy $q\le k-1$ and $q\le k-2$, respectively. Counting the normal derivatives in $S_aw_q$ and $D_aw_q$ shows that these sums involve normal derivatives of $b$ of order at most $s-1$. Thus the derivative $\partial_\nu^s b$ enters through $W_{k-1,\ell+1}$, and its coefficient is multiplied by $-1/2$ at each step.

When $k=1$, the term $-i^{-1}S_aw_0$ in $M_1$ contributes
$-i^{-1}(\partial_\nu^\ell a)\chi$. Solving
\eqref{eq:induc-k} therefore gives the coefficient
$-(2i\beta)^{-1}$ of $\partial_\nu^\ell a$.
For $k\ge2$, the direct terms involve normal derivatives of $a$ of order at most $s-2$. Hence the derivative
$\partial_\nu^{s-1}a$ also propagates through
$W_{k-1,\ell+1}$ with a factor of $-1/2$ at each step.
This establishes~\eqref{eq:Wkl}. Since $\beta>0$, the recursion uniquely determines all $W_{k,\ell}$.
\end{proof}

Using the determined boundary derivatives $W_{k,\ell}$, we construct the polynomial amplitudes of the following form
\begin{equation}
w_k(t,y)=\sum_{\ell=0}^{N-k}
\frac{y_1^\ell}{\ell!}W_{k,\ell}(t,y').
\label{eq:bd-finite-realization}
\end{equation}
Equations~\eqref{eq:W0l} and~\eqref{eq:Wkl} show that each
$W_{k,\ell}$ is smooth and supported in
$\operatorname{supp}\chi$. Hence its extension by zero
outside $U'$ is smooth.

\subsection{Remainder estimates and proof of Proposition~\ref{prop:boundary}}
Let the remainder \(r_\tau\) satisfy the following problem
\begin{equation*}
\begin{cases}
L_{a,b}^+ r_\tau = F_\tau & \text{in } Q, \\
r_\tau = 0 & \text{on } \Sigma, \\
r_\tau|_{t=0} = 0 & \text{in } \Omega,
\end{cases}
\end{equation*}
with source $F_\tau = -L_{a,b}^+ v_\tau^{\text{app}}.$ Set $w_{-1}=w_N=0$ and $M_0=0$. Expanding
$F_\tau=-L_{a,b}^+v_\tau^{\mathrm{app}}$ by
\eqref{eq:bd-conjugation} and~\eqref{eq:bd-history-expansion}
gives
\begin{equation*}
\begin{split}
F_\tau=e^{i\tau t+\tau\Phi}\vartheta\big(\tau^2bE_\Phi w_\tau^{\rm amp}
+\tau^2\mathcal K_\tau(aE_\Phi w_\tau^{\rm amp})+\sum_{k=0}^{N}\tau^{1-k}\mathcal F_k+\mathcal T_{N,N_{\rm time},\tau}
+\widetilde{\mathcal E}_{N_{\rm time},\tau}\big)+\mathcal C_{\vartheta,\tau},
\end{split}
\end{equation*}
where $\mathcal F_k=T_bw_k-P_b^+w_{k-1}-M_k$,
$\mathcal T_{N,N_{\rm time},\tau}$ collects the remaining algebraic terms from the finite temporal expansion, and
\begin{equation*}
\widetilde{\mathcal E}_{N_{\rm time},\tau}=\sum_{k=0}^{N-1}\left(\tau^{1-k}\mathcal E_{N_{\rm time},\tau}(S_aw_k)+\tau^{-k}\mathcal E_{N_{\rm time},\tau}(D_aw_k)\right).
\end{equation*}
Every term in $\mathcal T_{N,N_{\rm time},\tau}$ has a power at most $\tau^{-N}$. The cutoff commutator
$\mathcal C_{\vartheta,\tau}$ is supported where
$\nabla\vartheta\ne0$ and is exponentially small.

We use the normal decay estimate
\begin{equation}
\|\rho^qe^{-c\tau\rho}\|_{L^2(0,\varepsilon)} \le C_q\tau^{-q-1/2},\quad q\ge0.
\label{eq:bd-normal-decay}
\end{equation}
Since $E_\Phi=O(\rho^N)$ and the amplitudes are independent of $\tau$, \eqref{eq:bd-normal-decay} gives
\[
\|e^{\tau\Phi}\vartheta\,\tau^2bE_\Phi w_\tau^{\rm amp}\|_{L^2(Q)}=O(\tau^{3/2-N}).
\]
The amplitudes vanish near $t=0$. As $aE_\Phi$ is independent of time, integration by parts gives the identity
\begin{equation*}
\tau^2\mathcal K_\tau(aE_\Phi w_\tau^{\rm amp})
=-i\tau aE_\Phi w_\tau^{\rm amp}
+i\tau\mathcal K_\tau(aE_\Phi\partial_t w_\tau^{\rm amp}).
\end{equation*}
The bound for $\mathcal K_\tau$ on $L^2(Q)$ and
\eqref{eq:bd-normal-decay} therefore imply
\begin{equation*}
\bigl\|e^{\tau\Phi}\vartheta
\tau^2\mathcal K_\tau(aE_\Phi w_\tau^{\rm amp})\bigr\|_{L^2(Q)}
=O(\tau^{1/2-N}).
\end{equation*}

Equations~\eqref{eq:prin-amplitude} and \eqref{eq:bd-jet-equations} give $\mathcal F_k=O(\rho^{N-k})$ for $0\le k<N$. Moreover, the remaining coefficient $\mathcal F_N$ is bounded. Hence, for $0\le k\le N$, we have
\begin{equation*}
\bigl\|e^{\tau\Phi}\vartheta \tau^{1-k}\mathcal F_k\bigr\|_{L^2(Q)}=O(\tau^{1/2-N}).
\end{equation*}
The estimate \eqref{eq:bd-history-remainder} bounds the two temporal remainders for \(S_aw_k\) and \(D_aw_k\) by \(O(\tau^{-N_{\rm time}-1})\).
Since the time-independent weight \(\omega_\tau=e^{\tau\Phi}\vartheta\) commutes with the remainder operator, we apply \eqref{eq:bd-history-remainder} to \(\omega_\tau G\), with \(G=S_aw_k\) or \(G=D_aw_k\), and use \eqref{eq:bd-normal-decay} and the smoothness of \(G\) to bound the required time derivatives of \(\omega_\tau G\) by \(O(\tau^{-1/2})\) in \(L^2(Q)\).
Thus
\begin{align*}
\|e^{\tau\Phi}\vartheta\tau^{1-k}\mathcal E_{N_{\rm time},\tau}(S_aw_k)\|_{L^2(Q)}
&=O(\tau^{-N_{\rm time}-k-1/2}),\\
\|e^{\tau\Phi}\vartheta\tau^{-k}\mathcal E_{N_{\rm time},\tau}(D_aw_k)\|_{L^2(Q)}
&=O(\tau^{-N_{\rm time}-k-3/2}).
\end{align*}
Finally, since \(\operatorname{supp}(\nabla \vartheta)\) is disjoint from the boundary $\Gamma$, the commutator 
\(\mathcal{C}_{\vartheta,\tau}\) decays as \(O(\tau^{-L})\) for any fixed \(L > 0\).

Combining these estimates and choosing $N\ge M+4$ and $N_{\rm time}\ge N+3$ gives
\begin{align}
\|F_\tau\|_{L^2(Q)}\le C_{N,N_{\rm time}}\bigl(\tau^{3/2-N}+\tau^{1/2-N}+\tau^{-N_{\rm time}-1/2}+\tau^{-N_{\rm time}-3/2}\bigr)\le C_M\tau^{-M-2}.
\label{eq:bd-residual}
\end{align}
By Proposition~\ref{wellpose-Lab}, we obtain the estimate
\begin{equation}
\|r_\tau\|_{\mathbb W}+\|\partial_\nu r_\tau\|_{\mathbb N}\le C_M\tau^{-M-2}.
\label{eq:bd-correction}
\end{equation}
Thus \(V_\tau=v_\tau^{\mathrm{app}}+r_\tau\) solves \eqref{eq:velocity}. Its boundary velocity is \(h_\tau(t,y')=e^{i\tau(t+\psi(y'))}\chi(t,y')\).
Let \(\Sigma' := (0,T)\times U'\) denote the local lateral boundary in the coordinates \(y'\). The phase \(\psi\) is defined on \(U'\), and \(\chi\) is compactly supported in \(\Sigma'\). Since $\partial_\nu\Phi|_\Ga=1$ and the cutoff is identically 1 near the boundary, the expansion
\begin{equation*}
 e^{-i\tau(t+\psi)}\partial_\nu V_\tau =\tau\chi+\sum_{k=0}^M\tau^{-k}W_{k,1}+O_{L^2(\Sigma')}(\tau^{-M-1})
\end{equation*}
holds on \(\Sigma'\).
Removing the phase in \eqref{eq:velocityresponse} turns \(K_+\) into \(\mathcal K_\tau\), since \(\psi\) is independent of time. Applying \eqref{eq:bd-history-expansion} to the smooth terms and using the \(L^2\) boundedness of \(\mathcal K_\tau\) for the trace remainder gives
\begin{equation}
e^{-i\tau(t+\psi)}\mathcal R_{a,b}^+h_\tau =\tau\beta\chi+\sum_{m=0}^{M}\tau^{-m}B_m+O_{L^2(\Sigma')}(\tau^{-M-1}),
\label{eq:boundary-expansion}
\end{equation}
where
\begin{equation*}
B_m=\beta W_{m,1}+\alpha c_m\partial_t^m\chi
+\alpha\!\!\sum_{\substack{r,k\ge0\\r+k=m-1}}c_r\partial_t^rW_{k,1},
\qquad c_r=(-1)^ri^{-r-1}.
\end{equation*}
The sum is empty for \(m=0\). Substituting \eqref{eq:W0l} and \eqref{eq:Wkl} into this formula yields
\begin{equation}
B_m=\left(-\frac12\right)^{m+1}\partial_\nu^{m+1}b\,\chi+\left(-\frac12\right)^m i^{-1}\partial_\nu^ma\,\chi +H_m(J_\nu^mb,J_\nu^{m-1}a).
\label{eq:boundary-coef-Bm}
\end{equation}
Here $H_m$ depends only on normal derivatives of $b$ through order $m$ and of $a$ through order $m-1$, together with their tangential derivatives, the relevant time derivatives, and the fixed phase and boundary profile. In particular, the first two orders are given by
\begin{align}
B_0={}&-\tfrac12\partial_\nu b\,\chi-i\alpha\chi+H_0,
\label{eq:bd-first-flux}\\
H_0={}&-i\beta\langle d_\Ga\psi,d_\Ga\chi\rangle
-\tfrac i2\langle d_\Ga\beta,d_\Ga\psi\rangle\chi
-\tfrac\beta2(\Delta\Phi)|_\Ga\chi+\tfrac i2\chi,\notag\\
B_1={}&\tfrac14\partial_\nu^2b\,\chi+\tfrac i2\partial_\nu a\,\chi+H_1(J_\nu^1b,J_\nu^0a).
\label{eq:bd-second-flux}
\end{align}

We are now in a position to prove Proposition~\ref{prop:boundary}.

\begin{proof}[Proof of Proposition~\ref{prop:boundary}]
For $j=1,2$, we apply~\eqref{eq:boundary-expansion} to the
response $\mathcal R_{a_j,b_j}^+h_\tau$. Since
$f_\tau=K_+h_\tau\in\mathcal D_T$, the condition~\eqref{eq:dataequality}
gives $\mathcal R_{a_1,b_1}^+h_\tau
 =\mathcal R_{a_2,b_2}^+h_\tau$.
Dividing the difference of the expansions by $\tau$ and letting $\tau\to\infty$ yields
\[
(b_1-b_2)|_\Sigma\,\chi=0.
\]
Using the arbitrariness of the real-valued function $\chi$, we have $b_1=b_2$ on $\Sigma$.

We argue by induction. Suppose that \(J_\nu^m b_1=J_\nu^m b_2\) and \(J_\nu^{m-1}a_1=J_\nu^{m-1}a_2\) on \(\Sigma\).
When we subtract the two expansions through order \(\tau^{-m}\), the coefficients of all powers higher than \(\tau^{-m}\) and both \(H_m\) terms cancel. Multiplying by \(\tau^m\) and passing to the limit in
\(L^2(\Sigma')\) gives
\begin{equation*}
\left(-\frac12\right)^{m+1}
\partial_\nu^{m+1}(b_1-b_2)\,\chi+
\left(-\frac12\right)^m i^{-1}
\partial_\nu^m(a_1-a_2)\,\chi=0.
\end{equation*}
The coefficients and \(\chi\) are real-valued, so the real and imaginary parts give
\(\partial_\nu^{m+1}(b_1-b_2)\chi=0\) and
\(\partial_\nu^m(a_1-a_2)\chi=0\), respectively.
Since both the local coordinate neighborhood on \(\Gamma\) and the test function \(\chi\) are arbitrary, we obtain
\(\partial_\nu^{m+1}b_1=\partial_\nu^{m+1}b_2\) and
\(\partial_\nu^m a_1=\partial_\nu^m a_2\) on \(\Sigma\).
This completes the induction.

Differentiating these equalities along \(\Gamma\) and, for the
time-dependent coefficients \(b_j\), in time gives the corresponding identities for tangential and time derivatives. By smoothness, the identities for \(b_1\) and \(b_2\), including those for their derivatives, extend to \(t=0\) and \(t=T\). Thus, the proof is finished.
\end{proof}

\section{Construction of exponential solutions with memory}
\label{sec:exponential}
By Proposition~\ref{prop:boundary}, \(b_2-b_1\) vanishes to infinite order on \([0,T]\times\Gamma\). Extending \(b_1\) smoothly and \(b_2-b_1\) by zero gives smooth extensions \(\widetilde b_1\) and \(\widetilde b_2\) to a slightly larger smooth domain \(\widetilde\Omega\), with
\[
\widetilde b_1=\widetilde b_2
\quad\text{on }[0,T]\times(\widetilde\Omega\setminus\Omega).
\]
For a sufficiently small enlargement, both extensions are positive. By openness of simplicity and compactness of \([0,T]\), the extended metrics \(\widetilde b_j^{-1}(t,\cdot)e\) remain simple for every
\(t\in[0,T]\).

For each \(t\in[0,T]\), write
\begin{equation}
g_{b_j}(t)=b_j^{-1}(t,\cdot)e,\qquad
\widetilde g_j(t)=\widetilde b_j^{-1}(t,\cdot)e
\end{equation}
for the metrics on \(\overline\Omega\) and \(\overline{\widetilde\Omega}\), respectively.
We denote the distance between \(x,y\in\overline{\widetilde\Omega}\) for \(\widetilde g_j(t)\) by \(d_j^e(t;x,y)\), and the intrinsic distance between \(p,q\in\overline\Omega\) for \(g_{b_j}(t)\) by \(d_j^\Omega(t;p,q)\).

In this section, we aim to construct high-frequency exponential solutions for the following forward and backward equations,
\[
\begin{cases}
L_{a,b}^+ v_\tau = 0 \\
v_\tau(0) = 0
\end{cases}
\quad \text{and} \quad
\begin{cases}
L_{a,b}^- p_\tau = 0 \\
p_\tau(T) = 0
\end{cases}
\]
of the form
\[
v_\tau = e^{\Theta_\tau}(\mu_+ + R_{+,\tau}), \quad p_\tau = e^{-\Theta_\tau}(\mu_- + R_{-,\tau}),
\]
where $\tau$ is a large parameter, $\Theta_\tau = \tau^2 t + \tau \phi(t,x)$, and $R_{\pm,\tau}$ denote the remainder terms that tend to zero as $\tau \to \infty$ in suitable norms. The conjugated time kernel has $L^1$ norm $O(\tau^{-2})$, which compensates for the two spatial derivatives in the memory term.

\subsection{WKB expansion for phase and amplitudes}
Substituting \(e^{\pm\Theta_\tau}\mu_\pm\) into \(P_b^\pm\), we have
\begin{align}\label{eq:expan-P}
e^{\mp\Theta_\tau} P_b^\pm \big(e^{\pm\Theta_\tau} \mu_\pm\big) = \tau^2 \big(1 - b|\nabla\phi|^2\big)\mu_\pm \mp \tau \mathcal{T}_\phi^\pm \mu_\pm + P_b^\pm \mu_\pm,
\end{align}
where the operator $\mathcal{T}_\phi^\pm$ is defined by
\begin{equation*}
\begin{aligned}
    \mathcal{T}_\phi^\pm w = 2b\nabla \phi \cdot \nabla w + \big(\operatorname{div}(b\nabla \phi) \mp \phi_t\big) w.
\end{aligned}
\end{equation*}
To make the leading-order term in \eqref{eq:expan-P} vanish, we are required to solve the eikonal equation
\begin{equation}\label{eq:eikonal}
    b|\nabla\phi|^2 = 1 \quad\text{on } [0,T]\times\overline{\Omega}.
\end{equation}
For a single viscosity \(b\), write \(\widetilde b\) for its extension and \(\widetilde g(t)=\widetilde b^{-1}(t,\cdot)e\).
For \(x_0\in\partial\widetilde\Omega\), the distance phase function \(\phi(t,x)=d_{\widetilde g(t)}(x_0,x)\) solves \eqref{eq:eikonal}.
It follows that the amplitudes $\mu_{\pm}$ satisfy the transport equations
\begin{equation}\label{eq:transport}
    \mathcal{T}_\phi^\pm \mu_\pm = 0.
\end{equation}

For an inward-pointing Euclidean unit vector \(\omega\) at \(x_0\),
define the \(\widetilde g(t)\)-geodesic
\[
\gamma_{t,x_0,\omega}(r)
=\exp_{x_0}^{\widetilde g(t)}
 \bigl(r\sqrt{\widetilde b(t,x_0)}\,\omega\bigr),
\]
parametrized by arc length \(r\). In these coordinates, we write \(f(t,r,\omega)=f(t,\gamma_{t,x_0,\omega}(r))\). In particular, \(\phi_t\) denotes the time derivative at fixed \(x\), evaluated at \(x=\gamma_{t,x_0,\omega}(r)\).
For the extended coefficient \(\widetilde b\), the transport
equation \eqref{eq:transport} takes the form
\[
2\partial_r\mu_\pm+\bigl(\operatorname{div}(\widetilde b\nabla\phi)\mp\phi_t\bigr)\mu_\pm=0.
\]
Let \(V\) have compact closure in the set of inward directions. Choose \(r_*>0\) small enough that \(\gamma_{t,x_0,\omega}(r)\in\widetilde\Omega\setminus\overline\Omega\) for \(0<r\le r_*\), \(t\in[0,T]\), and \(\omega\in V\).
Choose \(\chi_\pm\in C_c^\infty((0,T)\times V)\) independently and prescribe 
\(\mu_\pm(t,r_*,\omega)=\chi_\pm(t,\omega)\).
Solving the transport equations in \(r\) gives
\[
\mu_\pm(t,r,\omega)=\chi_\pm(t,\omega)
\exp\!\left(-\frac12\int_{r_*}^{r}
\bigl(\operatorname{div}(\widetilde b\nabla\phi)\mp\phi_t\bigr)(t,s,\omega)\,\dd s\right).
\]
The compact angular support of \(\chi_\pm\) allows us to extend \(\mu_\pm\) smoothly by zero to directions outside \(V\).

\subsection{Remainder estimates}

Once the phase \(\phi\) and leading amplitudes \(\mu_\pm\) are constructed, the remainders \(R_{\pm,\tau} \in \mathbb{X}_\pm\) have to satisfy 
\begin{equation}\label{eq:AR=q}
A_\tau^\pm R_{\pm,\tau} = -F_{\pm,\tau},
\end{equation}
where 
\begin{equation*}
A_\tau^\pm = e^{\mp\Theta_\tau} L_{a,b}^\pm e^{\pm\Theta_\tau},\quad
F_{\pm,\tau} = P_b^\pm \mu_\pm - \mathcal{M}_\tau^\pm \mu_\pm,
\end{equation*}
with the conjugated memory operators
\begin{equation*}
\mathcal{M}_\tau^+ r(t) = e^{-\Theta_\tau(t)} \int_0^t D_a\big(e^{\Theta_\tau(s)} r(s)\big)\dd s, \quad
\mathcal{M}_\tau^- r(t) = e^{\Theta_\tau(t)} \int_t^T D_a\big(e^{-\Theta_\tau(s)} r(s)\big)\dd s.
\end{equation*}

To invert \(A_\tau^\pm\), we first study the conjugated parabolic operators
\begin{equation}\label{eq:def-B}
B_\tau^\pm= e^{\mp\Theta_\tau}P_b^\pm e^{\pm\Theta_\tau}
= P_b^\pm \mp \tau\mathcal T_\phi^\pm.
\end{equation}
We then treat the Volterra term in
\(A_\tau^\pm=B_\tau^\pm-\mathcal M_\tau^\pm\) as a perturbation.
We now state the following weighted energy estimate.

\begin{lemma}\label{prop:B,G}
Let \(\phi\) be the constructed phase function, satisfying
\eqref{eq:eikonal}. For each fixed \(\tau\), the operators
\(B_\tau^\pm\colon\mathbb X_\pm\to\mathbb Y\) are isomorphisms.
There exist constants \(C>0\) and \(\tau_0>0\) such that, for every
\(\tau\ge\tau_0\) and \(F\in\mathbb Y\), the unique solution
\(r\in\mathbb X_\pm\) of \(B_\tau^\pm r=F\) satisfies
\begin{equation}\label{eq:est-B}
\tau\|r\|_{\mathbb Y}
+\tau^{1/2}\|\nabla r\|_{\mathbb Y}
+\tau^{-1/2}\bigl(
  \|r\|_{\mathbb W}+\|\partial_\nu r\|_{\mathbb N}
\bigr)
\le C\|F\|_{\mathbb Y}.
\end{equation}
\end{lemma}

\begin{proof}
Since \(P_b^\pm:\mathbb X_\pm\to\mathbb Y\) are isomorphisms and the weights $e^{\pm \Theta_{\tau}}$ have bounded inverses, \(B_\tau^\pm\) are isomorphisms for each fixed $\tau$.

We consider first the forward operator \(B_\tau^+\), and  rewrite the equation \(B_\tau^+r = F\) as
\[
\partial_t r - \operatorname{div}(b\nabla r) - \tau \mathcal{T}_\phi^+ r = F.
\]
Take $\lambda>\|\phi_t\|_\infty+1$ and set $w=e^{\lambda\phi}$. Testing the equation by \(w\bar{r}\), taking the real part, and integrating by parts, we apply $r|_\Sigma=0$ and \eqref{eq:eikonal} to obtain
\begin{equation*}
 \operatorname{Re}\int_\Q w\overline r\,\mathcal T_\phi^+r\,\dd x\dd t
   =-\int_\Q(\lambda+\phi_t)w|r|^2\,\dd x\dd t.
\end{equation*}
\begin{align*}
 \operatorname{Re}\int_\Q w\overline r\,r_t
  &=\frac12\int_\Om w(T)|r(T)|^2 \dd x
       -\frac12\int_\Q w_t|r|^2\dd x\dd t,\\
 \operatorname{Re}\int_\Q w\overline r(-\operatorname{div}(b\nabla r))
  &=\int_\Q wb|\nabla r|^2\dd x\dd t
       -\frac12\int_\Q\operatorname{div}(b\nabla w)|r|^2\dd x\dd t.
\end{align*}
We sum these identities and drop the nonnegative endpoint term.
For sufficiently large \(\tau\), we absorb the lower-order terms, which are bounded by \(C\|r\|_{\mathbb Y}^2\), into the positive term
\(c\tau\|r\|_{\mathbb Y}^2\). Applying the Cauchy-Schwarz inequality to the source term then yields
\begin{equation*}
 c\tau\|r\|_{\mathbb Y}^2+c\|\nabla r\|_{\mathbb Y}^2
       \le C\|F\|_{\mathbb Y}\|r\|_{\mathbb Y}.
\end{equation*}
Hence, we can get
\begin{equation}\label{eq:est_r+Y}
\|r\|_{\mathbb Y}\le C\tau^{-1}\|F\|_{\mathbb Y}, \quad \|\nabla r\|_{\mathbb Y}\le C\tau^{-1/2}\|F\|_{\mathbb Y}.
\end{equation}
Since \(\mathcal T_\phi^+\) is a first-order operator with bounded coefficients, estimate~\eqref{eq:est_r+Y} gives
\[
\tau\|\mathcal T_\phi^+r\|_{\mathbb Y}
\le C\tau\bigl(\|\nabla r\|_{\mathbb Y}+\|r\|_{\mathbb Y}\bigr)\le C\tau^{1/2}\|F\|_{\mathbb Y}.
\]
Since $P_b^+r = F + \tau \mathcal{T}_\phi^+r$, the triangle inequality yields $\|P_b^+r\|_{\mathbb Y}\le C\tau^{1/2}\|F\|_{\mathbb Y}$. Applying the estimate \eqref{eq:est-z} with \(\lambda_0=1\), we have
\begin{equation}\label{eq:est_r+W}
 \|r\|_{\mathbb W}
 \le C\|(P_b^++1)r\|_{\mathbb Y}
 \le C\bigl(\|P_b^+r\|_{\mathbb Y}+\|r\|_{\mathbb Y}\bigr)
 \le C\tau^{1/2}\|F\|_{\mathbb Y}.
\end{equation}
Combining \eqref{eq:est_r+Y}, \eqref{eq:est_r+W}, and the estimate \eqref{eq:normaltrace} for normal trace yields \eqref{eq:est-B}.

For the backward problem, we test \(B_\tau^-r=F\) against
\(w\overline r\), with \(w=e^{-\lambda\phi}\). Integration by parts
in space gives
\[
\operatorname{Re}\int_Q w\overline r\,\mathcal T_\phi^-r\,\dd x\,\dd t=\int_Q(\lambda+\phi_t)w|r|^2\,\dd x\dd t.
\]
Since \(r(T,\cdot)=0\), integration by parts in time leaves only a nonnegative endpoint term at \(t=0\). The remaining estimates proceed as in the forward case and yield \eqref{eq:est-B}.
\end{proof}

Let \(G_\tau^\pm\) denote the inverses of \(B_\tau^\pm\). In the conjugation of \(D_a=\operatorname{div}(a\nabla)\) by \(e^{\pm\Theta_\tau}\),
the zeroth-, first-, and second-order coefficients have sizes
\(O(\tau^2)\), \(O(\tau)\), and \(O(1)\), respectively. To estimate \(\mathcal M_\tau^\pm\), we therefore equip
\(\mathbb Z_\tau=L^2(0,T;H^2(\Omega))\) with the scaled norm
\[
\|r\|_{\mathbb Z_\tau}=\|r\|_{\mathbb Y} +\tau^{-1}\|\nabla r\|_{\mathbb Y}
+\tau^{-2}\|D_x^2r\|_{\mathbb Y}.
\]

\begin{lemma}\label{lem:bound-M}
There exist constants \(C>0\) and \(\tau_0>0\) such that, for all \(\tau\ge\tau_0\),
\[
\|\mathcal M_\tau^\pm\|_{\mathbb Z_\tau\to\mathbb Y}\le C.
\]
Moreover, 
\[
\|\mathcal M_\tau^\pm G_\tau^\pm\|_{\mathbb Y\to\mathbb Y}
\le C\tau^{-1}.
\]
\end{lemma}

\begin{proof}
Let \(L = \|\phi_t\|_\infty\). For \(s \leq t\) and \(\tau \geq 2L + 1\), the forward Volterra kernel satisfies
\begin{equation*}
e^{-\Theta_\tau(t,x)+\Theta_\tau(s,x)} \leq e^{-\tau^2(t-s)/2},
\end{equation*}
which yields an \(L^1\)-norm bound of \(2\tau^{-2}\).
At each fixed time, conjugation of \(D_a\) gives
\begin{equation*}
e^{-\Theta_\tau(s)} D_a \big(e^{\Theta_\tau(s)} r\big)
= D_a r + \tau\big(2a\nabla\phi\cdot\nabla r + \operatorname{div}(a\nabla\phi)\,r\big) + \tau^2 a|\nabla\phi|^2 r.
\end{equation*}
Combining the conjugation formula with the kernel bound and applying Young's convolution inequality in time, we obtain
\[
\|\mathcal{M}_\tau^+ r\|_{\mathbb{Y}} \le C\tau^{-2} \big( \tau^2\|r\|_{\mathbb{Y}} + \tau\|\nabla r\|_{\mathbb{Y}} + \|D_a r\|_{\mathbb{Y}} \big) \le C\|r\|_{\mathbb{Z}_\tau}.
\]
The estimate for \(\mathcal{M}_\tau^-\) follows by the same argument for \(s \geq t\).

To estimate \(\mathcal M_\tau^\pm G_\tau^\pm\), we fix
\(F\in\mathbb Y\) and set \(r=G_\tau^\pm F\).
Applying \eqref{eq:est-B} to \(B_\tau^\pm r=F\) gives
\[
\|r\|_{\mathbb{Y}} \le C\tau^{-1}\|F\|_{\mathbb{Y}}, \quad \|\nabla r\|_{\mathbb{Y}} \le C\tau^{-1/2}\|F\|_{\mathbb{Y}}, \quad \|r\|_{\mathbb{W}} \le C\tau^{1/2}\|F\|_{\mathbb{Y}}.
\]
Using $\|D_x^2 r\|_{\mathbb{Y}} \le \|r\|_{\mathbb{W}}$, we obtain
\[
\|G_\tau^\pm F\|_{\mathbb{Z}_\tau} \le \|r\|_{\mathbb{Y}} + \tau^{-1}\|\nabla r\|_{\mathbb{Y}} + \tau^{-2}\|r\|_{\mathbb{W}} \le C\big(\tau^{-1} + \tau^{-3/2} + \tau^{-3/2}\big)\|F\|_{\mathbb{Y}} \le C\tau^{-1}\|F\|_{\mathbb{Y}}.
\]
Thus, we have $\|D_x^2 r\|_{\mathbb{Y}} \le C\|r\|_{\mathbb{W}}$. Combining this with the bound for $\mathcal{M}_\tau^\pm$ yields the desired composition estimate. 
\end{proof}

By Lemma \ref{lem:bound-M}, we can choose \(\tau\) sufficiently large such that  \(\|\mathcal{M}_\tau^\pm G_\tau^\pm\|_{\mathbb{Y} \to \mathbb{Y}} \le 1/2 < 1\). Then \(I - \mathcal{M}_\tau^\pm G_\tau^\pm\) is an isomorphism on \(\mathbb{Y}\). The factorization
\begin{equation*}
A_\tau^\pm = B_\tau^\pm - \mathcal{M}_\tau^\pm = (I - \mathcal{M}_\tau^\pm G_\tau^\pm) B_\tau^\pm,
\end{equation*}
implies that \(A_\tau^\pm: \mathbb{X}_\pm \to \mathbb{Y}\) are isomorphisms. Their inverses are given by
\begin{equation*}
\widetilde{G}_\tau^\pm = G_\tau^\pm \big(I - \mathcal{M}_\tau^\pm G_\tau^\pm\big)^{-1},
\end{equation*}
which satisfy the same asymptotic bounds as \(G_\tau^\pm\).

Based on the above discussion, we can obtain the following proposition regarding the construction of exponential solutions.
\begin{proposition}\label{prop:ex-solution}
Fix \(x_0\in\partial\widetilde\Omega\), and let \(\phi\) be the
associated distance phase satisfying \eqref{eq:eikonal}.
Let \(\mu_\pm\) be smooth solutions of \eqref{eq:transport} compactly supported in time in \((0,T)\), as constructed above.
For all sufficiently large \(\tau\), the equations
\(L_{a,b}^+v_\tau=0\) and \(L_{a,b}^-p_\tau=0\) in \(Q\) admit
solutions of the form
\[
v_\tau=e^{\Theta_\tau}(\mu_++R_{+,\tau}),\quad
p_\tau=e^{-\Theta_\tau}(\mu_-+R_{-,\tau}),
\]
with \(v_\tau(0,\cdot)=0\) and \(p_\tau(T,\cdot)=0\).
Here \(\Theta_\tau=\tau^2t+\tau\phi(t,x)\), and the remainders
\(R_{\pm,\tau}\in\mathbb X_\pm\) satisfy \eqref{eq:AR=q} and the estimate
\[
\tau\|R_{\pm,\tau}\|_{\mathbb Y}
+\tau^{1/2}\|\nabla R_{\pm,\tau}\|_{\mathbb Y}
+\tau^{-1/2}\bigl(\|R_{\pm,\tau}\|_{\mathbb W}
+\|\partial_\nu R_{\pm,\tau}\|_{\mathbb N}\bigr)\le C,
\]
where \(C>0\) is independent of \(\tau\).
\end{proposition}

\begin{proof}
Set \(R_{\pm,\tau}=-\tilde G_\tau^\pm F_{\pm,\tau}\). Then
\(R_{\pm,\tau}\in\mathbb X_\pm\) and \eqref{eq:AR=q} holds.
The amplitudes \(\mu_\pm\) are smooth and independent of \(\tau\), so
\[
\|P_b^\pm\mu_\pm\|_{\mathbb Y}+\|\mu_\pm\|_{\mathbb Z_\tau}\le C, \quad \tau\ge 1.
\]
Lemma~\ref{lem:bound-M} therefore gives
\(\|F_{\pm,\tau}\|_{\mathbb Y}\le C\) for all sufficiently large \(\tau\).
The bounds on \(\tilde G_\tau^\pm\) yield the stated remainder estimate.
The amplitudes vanish near \(t=0\) and \(t=T\), and
\(R_{+,\tau}(0,\cdot)=R_{-,\tau}(T,\cdot)=0\) by the definition of \(\mathbb X_\pm\).
Hence \(v_\tau(0,\cdot)=0\) and \(p_\tau(T,\cdot)=0\).
\end{proof}

Since $R_{\pm,\tau}|_\Sigma = 0$, the normal traces on $\Sigma$ reduce to
\begin{equation*}
\begin{aligned}
e^{-\Theta_\tau} \partial_\nu v_\tau &= \tau (\partial_\nu\phi)\mu_+ + \partial_\nu\mu_+ + \partial_\nu R_{+,\tau}, \\
e^{\Theta_\tau} \partial_\nu p_\tau &= -\tau (\partial_\nu\phi)\mu_- + \partial_\nu\mu_- + \partial_\nu R_{-,\tau}.
\end{aligned}
\end{equation*}
By Proposition~\ref{prop:ex-solution}, the remainder traces satisfy $\|\partial_\nu R_{\pm,\tau}\|_{\mathbb{N}} = O(\tau^{1/2})$.

\section{Proof of the main theorem}\label{sec:proof-theorem}
This section is devoted to proving the main theorem of this paper.

\subsection{Determination of the viscosity}\label{sec-deter-b}
Our first objective is to recover the metric $g_b(t)$ without assuming prior knowledge of the elastic coefficient in the interior. We begin by establishing an estimate for the Volterra memory terms.

\begin{lemma}\label{lem:history_estimate}
Let $\phi_j \in C^1([0, T] \times \partial\Omega)$ be real-valued and set $\Theta_j = \tau^2 t + \tau \phi_j$ for $j=1,2$.
Let $L = \max_j \|\partial_t \phi_j\|_{L^\infty(\Sigma)}$, $\kappa_\tau = \tau^2 - L\tau$.
Suppose that $c \ge 0$ and $E \subset \Sigma$ is a measurable set such that $\phi_2 - \phi_1 \ge c$ on $E$.
Then for all $\tau \ge 2L + 1$ and $q \in L^2(\Sigma)$,
\begin{align}
\|\mathbf{1}_E e^{-\Theta_2} K_+(e^{\Theta_1} q)\|_{L^2(\Sigma)} &\le e^{-c\tau} \kappa_\tau^{-1} \|q\|_{L^2(\Sigma)}, \label{eq:history_1} \\
\|\mathbf{1}_E e^{\Theta_1} K_-(e^{-\Theta_2} q)\|_{L^2(\Sigma)} &\le e^{-c\tau} \kappa_\tau^{-1} \|q\|_{L^2(\Sigma)}. \label{eq:history_2}
\end{align}
\end{lemma}

\begin{proof}
For $(t, x) \in E$ and $s \le t$, it follows from the definitions of $L$ and $E$ that
\begin{equation*}
\begin{aligned} 
-\Theta_2(t,x) + \Theta_1(s,x) &= \tau^2(s - t) + \tau\bigl(\phi_1(s, x) - \phi_1(t, x)\bigr) + \tau\bigl(\phi_1(t, x) -\phi_2(t, x)\bigr) \\
 &\le -\tau^2(t-s) + \tau L(t-s) - \tau c = -c\tau - \kappa_\tau(t - s). 
\end{aligned}
\end{equation*}
By Young's convolution inequality, we obtain
\begin{equation*}
\begin{aligned}
\|\mathbf{1}_E e^{-\Theta_2} K_+(e^{\Theta_1} q)\|_{L^2(\Sigma)} &\leq e^{-c\tau} \left\| \int_0^t e^{-\kappa_\tau(t-s)} |q(s,\cdot)| ds \right\|_{L^2(\Sigma)} \\
&\leq e^{-c\tau} \|e^{-\kappa_\tau(\cdot)}\|_{L^1(0,\infty)} \|q\|_{L^2(\Sigma)} \\
&= e^{-c\tau} \kappa_\tau^{-1} \|q\|_{L^2(\Sigma)}.
\end{aligned}
\end{equation*}
Similarly, \eqref{eq:history_2} follows by integrating backwards from \(t\) to \(T\).
\end{proof}

We now show that the map \(\Lambda_{a,b}^T\) uniquely determines \(b\).
\begin{proof}[Proof of $b_1=b_2$ in Theorem \ref{thm:main}]
We divide the proof into three steps.

\noindent \textbf{Step 1.} Fix $t_0 \in (0, T)$ and define
\begin{equation}\label{eq:max_dist}
\ell = \max_{\substack{x\in\partial\tilde{\Omega},\,  y\in\partial\Omega}} |d_2^e(t_0; x, y) - d_1^e(t_0; x, y)|.
\end{equation}
We argue by contradiction, assuming that $\ell > 0$.
Interchanging the two coefficient pairs if necessary, we can choose $(x_0, y_0)$ such that $d_2^e(t_0; x_0, y_0) - d_1^e(t_0; x_0, y_0) = \ell$, and define
\begin{equation*}
  \phi_j(t,x):=d^{\mathrm{e}}_j(t;x_0,x),\quad
  \Theta_j(t,x)=\tau^2t+\tau\phi_j(t,x),\quad j=1,2 .
\end{equation*}
The metrics \(\widetilde g_1(t)\) and \(\widetilde g_2(t)\) agree on \(\overline{\widetilde\Omega}\setminus\Omega\), and both are simple on \(\overline{\widetilde\Omega}\).
By the connected-intersection property of simple manifolds, every geodesic crossing \(\Omega\) enters and exits it exactly once, both transversely.

For $j=1,2$, let $\gamma_j$ be the minimizing geodesic from $x_0$ to $y_0$ in the metric $\widetilde g_j(t_0)$.
Parametrize \(\gamma_j\) by \(\widetilde g_j(t_0)\)-arc length, and let
\(\xi_j=\widetilde g_j(t_0)(\dot\gamma_j(0),\cdot)\) be its initial
covector. The first variation formula and the maximality in
\eqref{eq:max_dist} give
\[
0=\frac{\dd}{\dd s}\Big|_{s=0}
\bigl(d_2^e(t_0;x(s),y_0)-d_1^e(t_0;x(s),y_0)\bigr)
=(\xi_1-\xi_2)(w)
\]
for any smooth curve \(x(s)\) in \(\partial\widetilde\Omega\) with \(x(0)=x_0\) and \(\dot x(0)=w\). Thus \(\xi_1\) and \(\xi_2\) have the same restriction to 
\(T_{x_0}\partial\widetilde\Omega\).
With respect to the same metric at \(x_0\), both covectors have unit norm and positive inward normal components. Their tangential components therefore determine their normal components, so
\(\xi_1=\xi_2\) and \(\dot\gamma_1(0)=\dot\gamma_2(0)\).
Since the metrics agree on 
\(\widetilde\Omega\setminus\Omega\), uniqueness for the geodesic initial value problem implies that
\(\gamma_1\) and \(\gamma_2\) coincide up to their first intersection with \(\partial\Omega\). Denote this common intersection point by \(p\).

If \(p=y_0\), the two geodesics coincide, contradicting
\(\ell>0\). Since \(p\ne y_0\), uniqueness in the extended simple manifold identifies each segment from \(p\) to \(y_0\) with the corresponding intrinsic geodesic in \(\overline\Omega\). By strict convexity, these segments lie in \(\Omega\) between their endpoints and meet \(\Gamma\) transversely at \(p\) and \(y_0\).

By transversality at \(p\) and smooth dependence on \(t\) and \(x\), the geodesics from \(x_0\) to \(x\in\Gamma\) remain outside \(\Omega\) until they reach \(x\) for \((t,x)\) near \((t_0,p)\). The metrics agree along these geodesics, so their lengths and terminal covectors agree. The first variation formula then gives
\(\partial_\nu\phi_1=\partial_\nu\phi_2\), with a negative common value at \((t_0,p)\).
By continuity and the equality \(\phi_2(t_0,y_0)-\phi_1(t_0,y_0)=\ell\), we can choose an interval \(J\Subset(0,T)\) containing \(t_0\), disjoint open neighborhoods \(U_0,W_0\subset\Gamma\) of \(p,y_0\), respectively,
and \(c_0>0\) such that, for \(U=J\times U_0\) and \(W=J\times W_0\),
\begin{equation}\label{eq:phase_U}
\phi_1=\phi_2=:\phi,\quad
\partial_\nu\phi_1=\partial_\nu\phi_2=\partial_\nu\phi<-c_0\quad\text{on }U,
\end{equation}
\begin{equation}
\phi_2-\phi_1\ge\ell/2
\quad\text{on }W. \label{eq:phase_W}
\end{equation}

\noindent \textbf{Step 2.} Let \(\omega_j\) be the initial Euclidean unit direction of \(\gamma_j\). Choose smooth, nonnegative transport data for \(\mu_j\) that are
positive at \((t_0,\omega_j)\) and compactly supported in a small neighborhood of \((t_0,\omega_j)\). By transversality and smooth dependence of the intersections with \(\Gamma\), we have $\operatorname{supp}(\mu_j|_\Sigma)\subset U\cup W,$ $j=1,2$. Both amplitudes are positive at \((t_0,p)\), so \(\mu_1\mu_2>0\) on a nonempty open subset of \(U\).

Applying Proposition~\ref{prop:ex-solution} to the forward equation for \((a_1,b_1)\) and the backward equation for \((a_2,b_2)\), we can construct the following solutions
\[
v_1=e^{\Theta_1}(\mu_1+R_1),\quad
p_2=e^{-\Theta_2}(\mu_2+R_2),
\]
with \(R_1|_\Sigma=R_2|_\Sigma=0\).
Define \(q_1\) and \(q_2\) by
\begin{equation*}
q_1=\tau\partial_\nu\phi_1\mu_1+\partial_\nu\mu_1+\partial_\nu R_1,\quad
q_2=-\tau\partial_\nu\phi_2\mu_2+\partial_\nu\mu_2+\partial_\nu R_2.
\end{equation*} 
We have $\partial_\nu v_1 = e^{\Theta_1} q_1$ and $\partial_\nu p_2 = e^{-\Theta_2} q_2$, and $\|q_1\|_{L^2(\Sigma)} + \|q_2\|_{L^2(\Sigma)} \le C\tau$.
It follows from Proposition \ref{prop:boundary} that $\beta = b_1|_\Sigma = b_2|_\Sigma$ and $\alpha = a_1|_\Gamma = a_2|_\Gamma$.
Substituting $v_1, p_2$ into the boundary identity \eqref{eq:boundaryidentity}, we obtain
\begin{align}\label{eq:boundid-comp}
0  = \mathcal{B}_\tau := \mathcal{B}(v_1, p_2) =\mathcal{L}_\tau(\Sigma)+\mathcal{V}_{+,\tau}(\Sigma)
  -\mathcal{V}_{-,\tau}(\Sigma)
\end{align}
where
\begin{align*}
    \mathcal{L}_\tau(E) &= \int_E e^{\Theta_1-\Theta_2} \beta (\mu_2 q_1 - \mu_1 q_2) \dd S\dd t, \\
    \mathcal{V}_{+,\tau}(E) &= \int_{E} \alpha \mu_2 e^{-\Theta_2} K_+(e^{\Theta_1} q_1) \dd S\dd t, \\
    \mathcal{V}_{-,\tau}(E) &= \int_{E} \alpha \mu_1 e^{\Theta_1} K_-(e^{-\Theta_2} q_2) \dd S\dd t.
\end{align*}
Noting that the integrands vanish outside \(U \cup W\), we can split the identity into \(\mathcal{B}_\tau = \mathcal{B}_\tau(U) + \mathcal{B}_\tau(W)\), where \(\mathcal{B}_\tau(E) := \mathcal{L}_\tau(E) + \mathcal{V}_{+,\tau}(E) - \mathcal{V}_{-,\tau}(E)\).

On $U$, we have $\Theta_1 = \Theta_2$ by \eqref{eq:phase_U}.
Since $\partial_{\nu}\mu_j$ are bounded and $\|\partial_\nu R_j\|_{L^2(\Sigma)}=O(\tau^{1/2})$, we obtain
\begin{equation}\label{eq:LU_est}
    \mathcal{L}_\tau(U) = C_0 \tau + O(\tau^{1/2}), 
\end{equation}
where $ C_0 = 2 \int_U \beta \partial_\nu\phi \, \mu_1 \mu_2 \, dSdt < 0.$ Since \(\phi_1=\phi_2\) on \(U\), Lemma~\ref{lem:history_estimate}
with \(E=U\) and \(c=0\) gives
\[
|\mathcal V_{+,\tau}(U)|+|\mathcal V_{-,\tau}(U)|
\le C\tau^{-1}.
\]
Analogously on $W$, \eqref{eq:phase_W} and Lemma \ref{lem:history_estimate} yield
\begin{equation*}
  |\mathcal{B}_\tau(W)|\le C\tau e^{-\ell\tau/2}.
\end{equation*}
Combining these estimates, we arrive at the identity
\begin{equation}\label{eq:boundary-asymptotic}
    0 = C_0 \tau + O(\tau^{1/2}) + O(\tau e^{-\ell\tau/2}).
\end{equation}
For sufficiently large \(\tau\), the right-hand side of
\eqref{eq:boundary-asymptotic} is strictly negative, a contradiction. Hence \(\ell=0\), and
\begin{equation}\label{eq:dist_equal_e}
d_1^e(t_0;x,y)=d_2^e(t_0;x,y)
\quad\text{for all }(x,y)\in
\partial\tilde{\Omega}\times\partial\Omega.
\end{equation}

\noindent \textbf{Step 3.} Let $y, z \in \partial\Omega$ with $y \neq z$, and extend the intrinsic $g_{b_1}(t_0)$-geodesic from $y$ to $z$ backwards until it meets $\partial\tilde\Omega$ at a point $x$.
By the simplicity of \((\overline{\widetilde\Omega},\widetilde g_1(t_0))\), the extended geodesic from \(x\) to \(z\) and its subsegments are minimizing.
The portion from \(y\) to \(z\) has length \(d_1^\Omega(t_0;y,z)\) by construction, so
\begin{equation*}
d_1^e(t_0; x, z) = d_1^e(t_0; x, y) + d_1^\Omega(t_0; y, z).
\end{equation*}
On the other hand, the triangle inequality for $\widetilde g_2(t_0)$ gives
\begin{equation*}
d_2^e(t_0; x, z) \le d_2^e(t_0; x, y) + d_2^\Omega(t_0; y, z).
\end{equation*}
By \eqref{eq:dist_equal_e}, we obtain $d_1^\Omega(t_0; y, z) \le d_2^\Omega(t_0; y, z)$.
Interchanging the roles of \(g_{b_1}(t_0)\) and \(g_{b_2}(t_0)\) gives the reverse inequality
\(d_2^\Omega(t_0;y,z)\le d_1^\Omega(t_0;y,z)\).
Thus the intrinsic boundary distances agree:
\[
d_1^\Omega(t_0;y,z)=d_2^\Omega(t_0;y,z)
\quad\text{for all }(y,z)\in\partial\Omega\times\partial\Omega.
\]
By the boundary rigidity theorem of simple conformally Euclidean manifolds \cite[Theorem~11.4.1]{PaternainSaloUhlmann2023}, we obtain
\begin{equation*}
b_1(t_0, x) = b_2(t_0, x) \quad \forall x \in \bar{\Omega}.
\end{equation*}
Since $t_0\in(0,T)$ is arbitrary and $b_j$ are continuous, it follows that $b_1=b_2$ on $[0,T]\times \bar{\Omega}$.
\end{proof}

\subsection{Determination of the elastic coefficient}\label{sec:elasticity} 
To complete the proof of Theorem \ref{thm:main}, it remains to extract the elastic coefficient from the memory term in the mixed identity \eqref{eq:mixedidentity}.
We will assume throughout this section that the conditions of Theorem \ref{thm:main} are satisfied and note in view of the last section that it holds:
\begin{equation*}
b(t,x) := b_1(t,x) = b_2(t,x), \quad \text{in } Q.
\end{equation*}
We use the integral identity \eqref{eq:mixedidentity} to show that,
for each \(t\in(0,T)\), the scalar geodesic ray transform of \((a_1-a_2)/b(t,\cdot)\) with respect to \(g_b(t)\) vanishes.

To control the \(\tau\)-dependent remainders, we will show that the rescaled Volterra operators are uniformly bounded on \(L^2(Q)\) for large \(\tau\) and converge strongly to the identity as
\(\tau\to\infty\).

\begin{lemma}
\label{lem:approxidentity}
Let $\phi\in C^1([0,T]\times\overline\Om)$ be real-valued and set $L=\|\phi_t\|_{L^\infty(\Q)}$. For $\tau>L$, define
\begin{equation}
 H_\tau F(t,x)=\tau^2\int_0^t
 e^{-\tau^2(t-s)+\tau(\phi(s,x)-\phi(t,x))}F(s,x) \dd s.
 \label{eq:approxidentity}
\end{equation}
Then
\begin{equation}
 \|H_\tau\|_{L^2(\Q)\to L^2(\Q)}
 \le \frac{\tau}{\tau-L}
 \label{eq:approxidentitylimit}
\end{equation}
and $H_\tau F\to F$ in $L^2(\Q)$ for every $F\in L^2(\Q)$.
The same assertions hold componentwise for vector-valued functions. In particular, if $F_\tau\to F$ in $L^2(\Q)$, then $H_\tau F_\tau\to F$ in $L^2(\Q)$.
\end{lemma}

\begin{proof}
For $s\le t$, we have $\phi(s,x)-\phi(t,x) \le L(t-s).$ Hence, for almost every fixed $x\in \Omega $, Young's convolution inequality in the time variable yields
\begin{equation*}
     \|H_\tau F(\cdot,x)\|_{L^2(0,T)} \le \int_0^\infty \tau^2 e^{-(\tau^2-\tau L)r}\dd r \|F(\cdot,x)\|_{L^2(0,T)} = \frac{\tau}{\tau-L} \|F(\cdot,x)\|_{L^2(0,T)}.
\end{equation*}
Integrating the square of both sides over $x \in \Omega$ and applying Fubini's theorem gives the global bound\eqref{eq:approxidentitylimit}. Next, we define the integral operator
\begin{equation*}
 J_\tau F(t,x)=\tau^2\int_0^t e^{-\tau^2(t-s)}F(s,x) \dd s.
\end{equation*}
By extending $F$ by zero outside $(0,T)$, $J_\tau$ is identified as a classical one-sided approximate identity.
Thus, we have the strong convergence $J_\tau F\to F$ in $L^2(\Q)$.
A similar convolution argument yields
\begin{equation}
 \|H_\tau-J_\tau\|_{L^2\to L^2}
 \le\int_0^\infty\tau^2e^{-\tau^2r}(e^{\tau Lr}-1) \dd r
 =\frac{L}{\tau-L}.
 \label{eq:approxidentitycomparison}
\end{equation}
Applying the triangle inequality, we obtain
\begin{equation*}
\| H_\tau F - F \|_{L^2} \le \| H_\tau - J_\tau \|_{L^2\to L^2} \| F \|_{L^2} + \| J_\tau F - F \|_{L^2} \to 0.
\end{equation*}
Finally, by the uniform bound \eqref{eq:approxidentitylimit}, this strong convergence extends to parameter-dependent inputs $F_{\tau} \to F$.
\end{proof}

At the end of this section, we are now in a position to complete the proof of Theorem~\ref{thm:main}.
\begin{proof}[Proof of $a_1=a_2$ in Theorem \ref{thm:main}]
Write $\delta a=a_{1}-a_{2}$.
Fix an arbitrary $t_0\in(0,T)$ and a nontrivial maximal $g_b(t_0)$-geodesic \(\gamma\) in $\overline{\Omega}$.
Since the extensions of \(b_1\) and \(b_2\) also agree outside \(\Omega\), write \(\widetilde b\) for their common value and \(\widetilde g(t)=\widetilde b^{-1}(t,\cdot)e\) for the extended metric.
Extend \(\gamma\) backward in the common extended metric \(\widetilde g(t_0)\) to a point \(x_0\in\partial\widetilde\Om\).

By Proposition \ref{prop:ex-solution}, for this exterior pole $x_0$ and the associated distance phase function $\phi(t,x)$ satisfying $b|\nabla\phi|^2=1$, we can construct exponential solutions $v_1$ and $p_2$ solving $L_{a_1,b}^+ v_1=0$ and $L_{a_2,b}^-p_2 = 0$ respectively, which take the form
\begin{equation}
 v_1=e^{\Theta_\tau}(\mu_++R_1),\quad
 p_2=e^{-\Theta_\tau}(\mu_-+R_2),\quad
 \Theta_\tau=\tau^2t+\tau\phi.
 \label{eq:elasticgo}
\end{equation}
Substituting these solutions into the mixed volume identity \eqref{eq:mixedidentity} in Lemma \ref{lem:identities} gives
\begin{equation}
 \int_\Q\delta a\,\nabla(K_+v_1)\cdot\nabla p_2 \dd x\dd t=0.
 \label{eq:elasticidentity}
\end{equation}
Direct computation yields
\begin{align*}
     \tau e^{-\Theta_\tau}\nabla(K_+v_1)
&=H_\tau\left(\nabla\phi(\mu_++R_1)+\tau^{-1}\nabla(\mu_++R_1)\right),\\
\tau^{-1}e^{\Theta_\tau}\nabla p_2
 &=-\nabla\phi(\mu_-+R_2)+\tau^{-1}\nabla(\mu_-+R_2).
\end{align*}
By the remainder estimates in Proposition~\ref{prop:ex-solution}
and Lemma~\ref{lem:approxidentity}, we have
\[
\tau e^{-\Theta_\tau}\nabla(K_+v_1)\to\mu_+\nabla\phi,
\quad \tau^{-1}e^{\Theta_\tau}\nabla p_2\to-\mu_-\nabla\phi
\]
strongly in \(L^2(Q)\) as \(\tau\to\infty\). Taking the limit as $\tau \to \infty$ in \eqref{eq:elasticidentity} and using the eikonal equation, we obtain
\begin{equation}
 \int_\Q\frac{\delta a}{b}\,\mu_+\mu_- \dd x\dd t=0.
 \label{eq:elasticweighted}
\end{equation}

It remains to remove the transport weight. Let $(r,\omega)$ be the geodesic polar coordinates centered at $x_0$, with the volume element $dx = J_e(t,r,\omega)\,\dd r \dd\omega$. We note that along the geodesic $\gamma_{t,x_0,\omega}$, the differential operators can be expressed as
\begin{equation*}
\widetilde b\nabla\phi\cdot\nabla = \partial_r,\quad
\operatorname{div}_e(\widetilde b\nabla\phi) = \partial_r \log J_e.
\end{equation*}
Then the transport equations \eqref{eq:transport} can be rewritten as
\begin{equation*}
2\partial_r \mu_\pm + \bigl(\partial_r \log J_e \mp \phi_t\bigr)\mu_\pm = 0.
\end{equation*}
Multiplying the transport equations for \(\mu_+\) and \(\mu_-\) by \(\mu_-\) and \(\mu_+\), respectively, and adding the resulting equations gives
\begin{equation*}
\partial_r(\mu_+\mu_-) + (\partial_r \log J_e)\mu_+\mu_- = 0,
\end{equation*}
which implies
\begin{equation}
\partial_r\bigl(J_e \mu_+\mu_-\bigr) = 0.
\label{eq:transportproduct}
\end{equation}  

Let \(I\Subset(0,T)\) be an interval containing \(t_0\), and let \(V\) be a small neighborhood of the initial direction of the extended geodesic at \(x_0\). Choose \(r_*>0\) so that the initial segment \(0<r\le r_*\) lies in 
\(\widetilde\Omega\setminus\overline\Omega\) for all \((t,\omega)\in I\times V\). After shrinking \(I\) and \(V\), let \(r_-(t,\omega)\) and \(r_+(t,\omega)\) be the entry and exit arc-length parameters of the ray in \(\Omega\). The crossings are transverse, so these parameters depend smoothly on \((t,\omega)\) by the implicit function theorem.
For arbitrary
\(\chi_\pm\in C_c^\infty(I\times V)\), set
\[
\mu_\pm(t,r_*,\omega)
=J_e(t,r_*,\omega)^{-1/2}\chi_\pm(t,\omega).
\]
Thus \(J_e\mu_+\mu_-=\chi_+\chi_-\) along the ray from \(r_*\) through \(\Omega\). Substituting the product into \eqref{eq:elasticweighted} gives
\begin{equation}
\int_{I\times V}\chi_+\chi_- \Big( \int_{r_-(t,\omega)}^{r_+(t,\omega)} \frac{\delta a(\gamma_{t,x_0,\omega}(r))}{b(t,\gamma_{t,x_0,\omega}(r))}dr \Big)\dd \omega \dd t = 0.
\label{eq:elasticraytests}
\end{equation}
For any \(\chi_+\in C_c^\infty(I\times V)\), choose
\(\chi_-\in C_c^\infty(I\times V)\) equal to one on a neighborhood of \(\operatorname{supp}\chi_+\). Equation~\eqref{eq:elasticraytests}
then implies that the inner integral vanishes as a distribution on \(I\times V\). By continuity in \((t,\omega)\), it vanishes pointwise. In particular, the integral along \(\gamma\) at time \(t_0\) is zero. Since \(\gamma\) was an arbitrary nonconstant maximal \(g_b(t_0)\)-geodesic, we obtain
\[
I_{g_b(t_0)}\!\left(\frac{\delta a}{b(t_0,\cdot)}\right)=0.
\]
The manifold \((\overline{\Omega},g_b(t_0))\) is simple and has dimension \(n\ge2\). Injectivity of the scalar geodesic ray transform~\cite[Theorem~4.7.7]{PaternainSaloUhlmann2023} therefore gives \(\delta a/b(t_0,\cdot)=0\). Since \(b(t_0,\cdot)>0\), we conclude that \(\delta a=0\).
Therefore the proof of Theorem~\ref{thm:main} is complete.
\end{proof}

\section{Numerical experiments}\label{sec:numerics}

In this section, we present some numerical experiments of reconstructing the elastic coefficient \(a(x)\) and the time-dependent viscosity \(b(t,x)\) by time-domain full-waveform inversion. For each choice of boundary displacement, we solve the Kelvin--Voigt equation and compute moments of the resulting total traction. This gives a reduced optimization problem in the material coefficients. We refer
to~\cite{Plessix2006} for the adjoint-state approach and
to~\cite{PakravanKangNewtson2016} for related viscoelastic waveform inversion.

\subsection{Traction data and reduced objective}

Let \(h_\ell\) be a boundary velocity, \(f_\ell=K_+h_\ell\) the associated displacement, and \(k_m\) a boundary test function. The corresponding traction moment is given by
\begin{equation}\label{eq:num-moments}
 F_{\ell m}(a,b)
 =\int_\Sigma(\Lambda_{a,b}^T f_\ell)k_m\,\dd S\dd t.
\end{equation}
The boundary input determines the traction response, and the test function determines how we weight that response in space and time. For some inputs, we gradually change the boundary displacement and then keep it fixed. We include moments of the traction measured while the displacement is fixed in the fit.

We use the positive parametrization
\begin{equation}\label{eq:num-parameters}
 a_\theta=a_{\rm ref}e^{A_\theta},\quad
 b_\theta=b_{\rm ref}e^{B_\theta},\quad
 A_\theta=\sum_{i=1}^{265}\theta_i^a\psi_i,\quad
 B_\theta=\sum_{j=1}^{48}\theta_j^b\Psi_j .
\end{equation}
Here \(a_{\rm ref}=1\) and \(b_{\rm ref}=0.4\). The functions \(\psi_i\) are local tensor-product cubic B-splines restricted to the disk. The viscosity basis consists of the six spatial polynomials
\(1,x,y,x^2-y^2,2xy,x^2+y^2-\tfrac12\) multiplied by the first eight Legendre polynomials in \(2t/T-1\). This extends the temporal expansion from degree four to degree seven.

Let \(F_{\mathrm i}\) and \(F_{\mathrm d}\) be the moment maps on the inversion and data grids. For the full spline fit, we initialize at the coarse-grid estimate \(\theta_{\rm p}\) and evaluate the grid correction there.
We use a regularization penalty \(\mathcal R(\theta)\) on the coefficient parameters, specified in the next subsection, with weight \(\lambda>0\). Given the observations \(Y^\delta\), the residual and the reduced objective are
\begin{equation}\label{eq:num-objective}
 \begin{aligned}
 c&=F_{\mathrm d}(\theta_{\rm p})-F_{\mathrm i}(\theta_{\rm p}),\quad r(\theta)=W\operatorname{vec}
       \bigl(F_{\mathrm i}(\theta)+c-Y^\delta\bigr),\\
 \mathcal J_\lambda(\theta)
   &=\tfrac12\|r(\theta)\|_2^2+\lambda\mathcal R(\theta).
 \end{aligned}
\end{equation}
We keep the correction \(c\) fixed throughout each continuation run. For a group of \(L_g\) inputs, we multiply row \(\ell\) of the residual by 
\(\sqrt{1+\delta^2}/(\sqrt{2L_g}\,\|Y^\delta_{\ell,:}\|_2)\).
These weights account for the size of each observed row and give equal weight in the fit to the two groups: observations over the full time window and observations taken while the boundary displacement is fixed. Let 
\(v_\ell=\partial_tu_\ell\) and let \(p_m\) solve the terminal adjoint equation \(L^-_{a,b}p_m=0\) with boundary value \(k_m\) and zero terminal data. Applying~\eqref{eq:green} to the linearized state gives
\begin{equation}\label{eq:num-derivative}
 DF_{\ell m}(a,b)[\delta a,\delta b]
 =\int_Q(\delta a\nabla u_\ell+\delta b\nabla v_\ell)
         \cdot\nabla p_m\,\dd x\dd t .
\end{equation}
The explicit variation of the traction cancels the boundary term from the linearized divergence source. For parameter variations we substitute
\(\delta a=a_\theta\delta A_\theta\) and
\(\delta b=b_\theta\delta B_\theta\) in~\eqref{eq:num-derivative}.

\subsection{Joint updates and choice of regularization parameter}

In the full spline fit, \(\mathcal R\) includes a quadratic amplitude term and a pseudo-Huber gradient penalty on \(A_\theta\), a squared gradient penalty on \(a_\theta\), and a weighted squared Hessian penalty on \(A_\theta\). The Hessian weight is
\((1+|\nabla A_{\theta_{\rm p}}|^2/3^2)^{-1}\), which allows steeper transitions near the contours of the initial fit. The viscosity penalty is quadratic in \(B_\theta\), its spatial gradient, and its time derivative. Restricted elastic updates use directions derived from smoothed coarse-grid estimates. These directions describe contour shifts, contrast changes, and low-order background variations. Each restricted fit retains all 48 viscosity parameters and uses a quadratic penalty in the normalized reduced coordinates.

At each iterate \(\theta_n\), we use the discrete adjoint to compute the gradient \(g_n\) of the reduced objective. The damped direction solves a quadratic subproblem
with objective
\begin{equation}\label{eq:num-step}
 g_n^T d+\tfrac12 d^T
 \bigl(H_n+\lambda R_n+\mu_n D_n\bigr)d .
\end{equation}
Here \(H_n\) is the Gauss--Newton data curvature or its secant update, \(R_n\) is a positive approximation of the penalty curvature, and \(D_n\) is a positive diagonal scaling. We retain the cross block between the elastic and viscosity parameters and solve for both increments at the same step. The full fit constrains \(\theta_n+d\) to \([-1.2,1.2]^{313}\). The restricted fit uses directions in the reduced space and accepts trial coefficients in \([-1.6,1.6]^{313}\). We adjust the damping parameter according to the ratio of the actual decrease in the objective to the decrease predicted by the local model. 

We solve the inverse problem for selected values of
\(\lambda_j=1.5\times10^{-4}3^{-j}\) in decreasing order, using each fit to initialize the next. The discrepancy test accounts for noise and a grid allowance estimated from \(c\) in \eqref{eq:num-objective}. We compare stationary fits satisfying the discrepancy test on reserved waveforms from both groups. For noiseless data, a weighted difference between predictions on two grids sets the comparison threshold. For noisy data, we favor stronger regularization when score differences are within the estimated noise variation.

\subsection{Discretization and experimental setting}

We take \(\Omega\) to be the unit disk and \(T=3\). The spatial meshes triangulate the same inscribed regular 64-gon and use continuous piecewise quadratic finite elements. The data grid has 3329 nodes and 384 time steps. Five selected fits use the full spline space on a separate grid with 1777 nodes and 288 time steps. The restricted fit for the annulus at \(2\%\) noise uses 39 elastic update directions, 865 nodes, and 224 time steps.

We apply the implicit midpoint rule to the displacement and velocity. Let \(\mathsf M\), \(\mathsf A\), and \(\mathsf B_{n-1/2}\) be the mass and stiffness matrices for unit density, \(a_\theta\), and \(b_\theta(t_{n-1/2},\cdot)\). The displacement increment is \(\Delta t\,\mathbf v^{n-1/2}\), and the nodal reaction is
\begin{equation}\label{eq:num-state}
 \mathbf T^{n-1/2}=
 \frac{\mathsf M(\mathbf v^n-\mathbf v^{n-1})}{\Delta t}
 +\mathsf B_{n-1/2}\mathbf v^{n-1/2}
 +\mathsf A\mathbf u^{n-1/2}.
\end{equation}
The superscript \(n-1/2\) denotes the average of the two endpoint values for the states. The reaction vanishes at the interior nodes. We approximate \eqref{eq:num-moments} by pairing the boundary residual with the nodal test values and summing in time. Linearizing the time-stepping and traction equations gives the discrete tangent, whose transpose is the discrete adjoint.

The full-window group consists of 55 boundary inputs with 336 moments per input. The tests pair 21 real Fourier modes of orders 0 through 10 with 16 smooth time windows centered at equally spaced times from 0.075 to 2.925. The additional group consists of 20 displacement ramps. For each ramp input, we retain \(126\) traction moments from the last
six time windows, all within the interval of fixed boundary displacement. For validation, we reserve seven inputs with observations over the full time window and sixteen with observations taken while the displacement is fixed. For each observed row of clean data \(Y_{\ell,:}\), we draw an independent standard Gaussian vector \(\xi_\ell\) of the corresponding length and set
\begin{equation}\label{eq:num-noise}
Y^\delta_{\ell,:}=Y_{\ell,:}
+\delta\|Y_{\ell,:}\|_2\frac{\xi_\ell}{\|\xi_\ell\|_2},
\quad \delta\in\{0,0.02\}.
\end{equation}
We use one noise realization per case at \(2\%\).

\subsection{Reconstruction results}

The three elastic profiles contain two elliptic inclusions, a rotated kite, and a five-lobed annulus, respectively. The corresponding viscosity fields vary in amplitude, orientation, and location, respectively. All target fields and all six final reconstructions are strictly positive on their respective domains. The spline partition of unity bounds
\(A_\theta\) between the smallest and largest elastic coordinates. The spatial factors and the Legendre polynomials  are uniformly bounded on their domains, which bounds \(B_\theta\). These bounds and \eqref{eq:num-parameters} give positive lower bounds for the two reconstructed fields.

\begin{figure}[p]
 \centering
 \includegraphics[width=.80\textwidth]{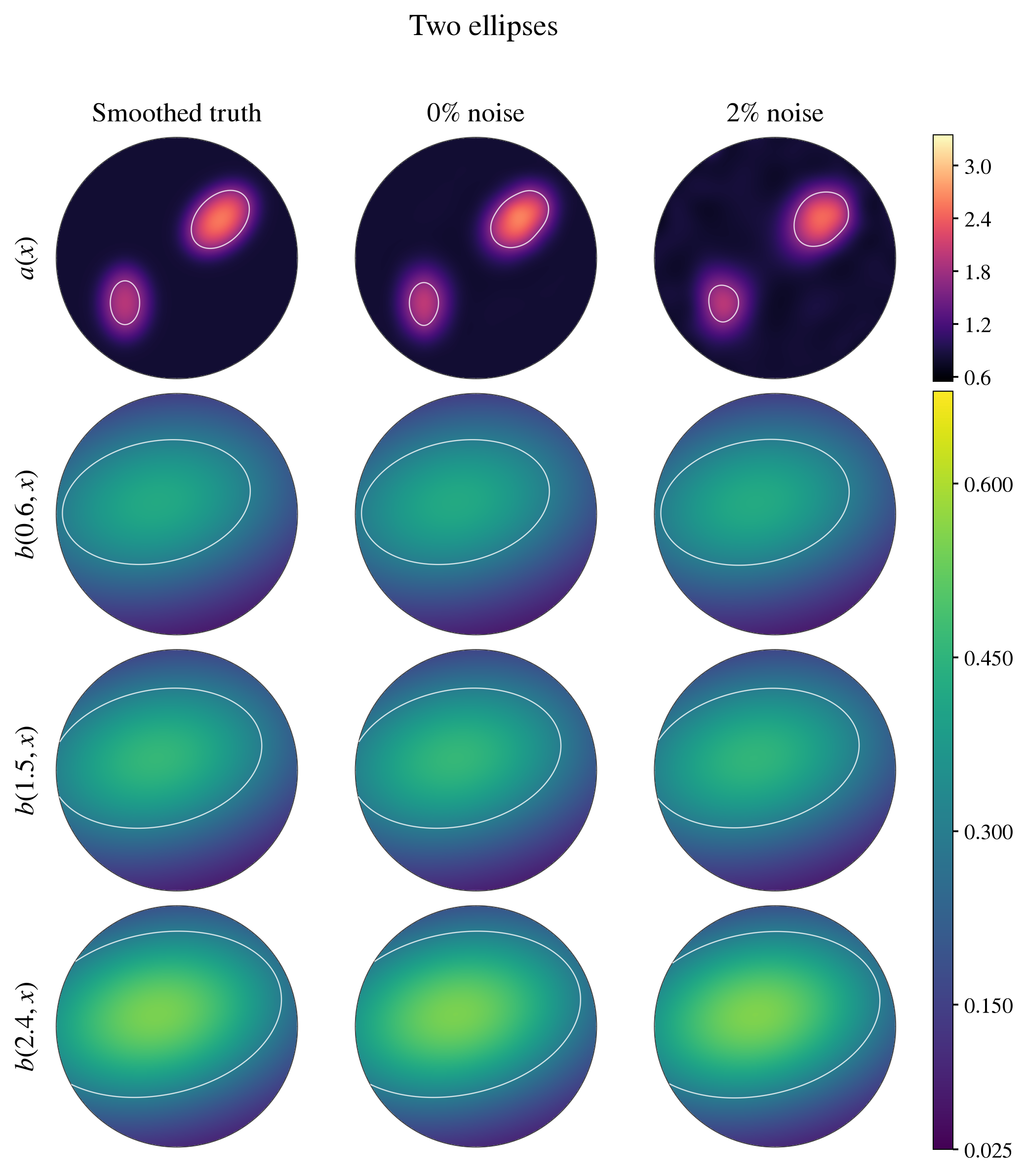}
 \caption{Case 1. The elastic field has two elliptic inclusions. The
 columns show the truth and the reconstructions from noiseless
 and \(2\%\) noisy data. The upper row shows \(a(x)\) and the remaining
 rows show \(b(t,x)\) at \(t=0.6,1.5,2.4\). Each field uses a common
 color scale. The contours mark \(a=1.65\) and \(b=0.30\).}
 \label{fig:num-ellipses}
\end{figure}

\begin{figure}[p]
 \centering
 \includegraphics[width=.88\textwidth]{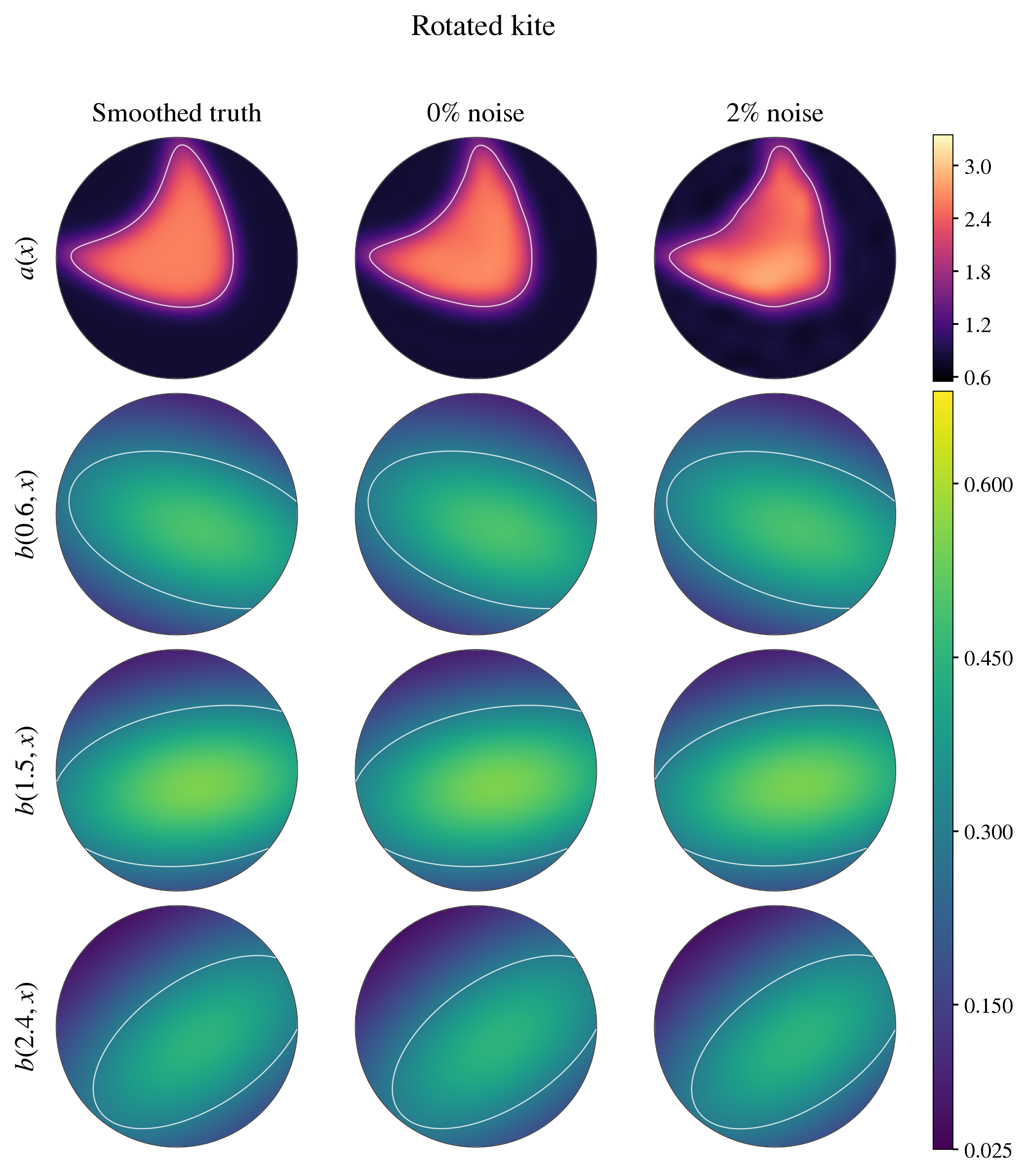}
 \caption{Case 2. Reconstruction of the rotated kite and the viscosity
 with a rotating principal axis. The columns, time slices, and color
 scales follow Figure~\ref{fig:num-ellipses}.}
 \label{fig:num-kite}
\end{figure}

\begin{figure}[p]
 \centering
 \includegraphics[width=.80\textwidth]{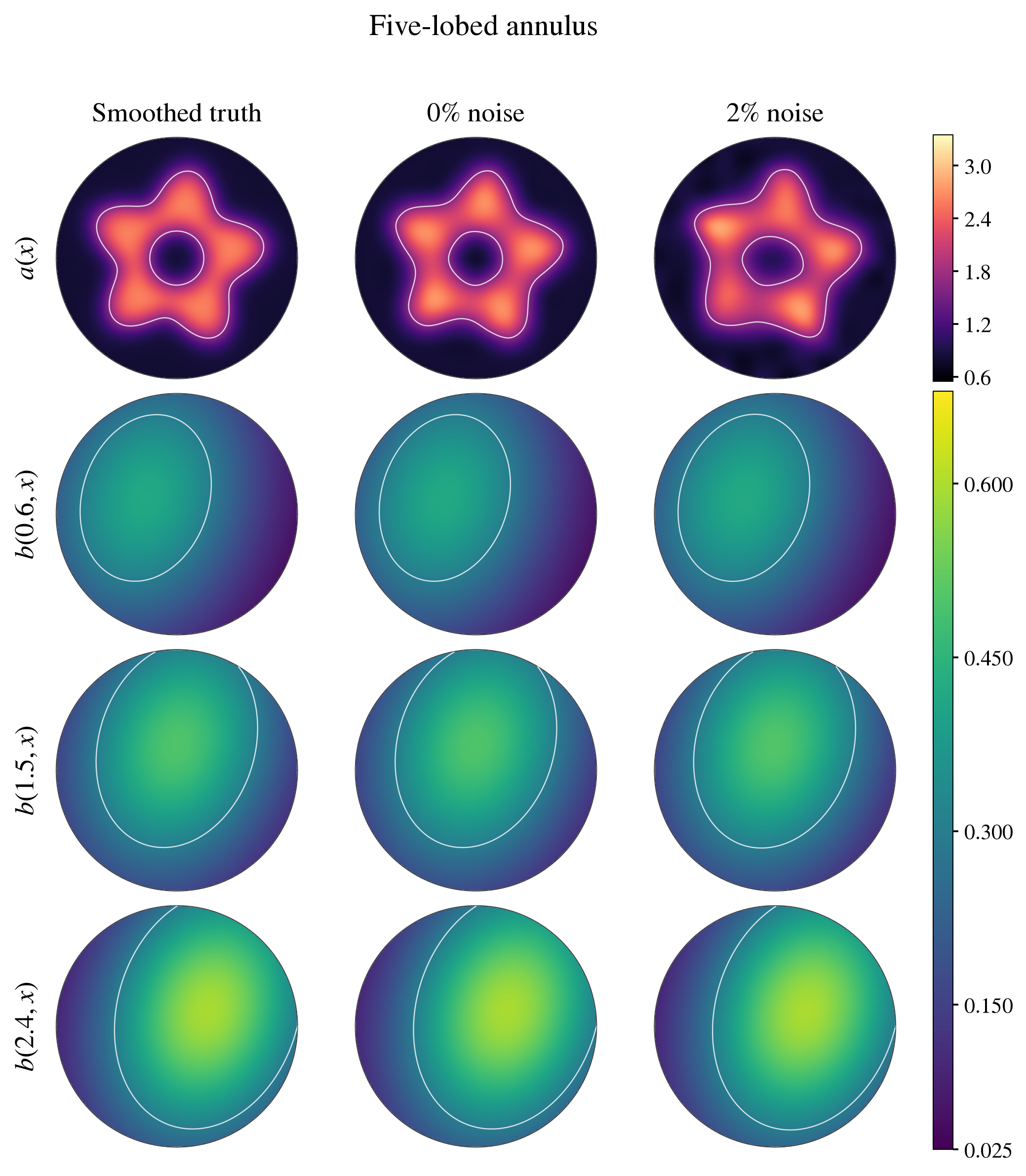}
 \caption{Case 3. Reconstruction of the five-lobed annulus and the
 viscosity with a moving center. The columns, time slices, and color
 scales follow Figure~\ref{fig:num-ellipses}.}
 \label{fig:num-flower}
\end{figure}

Figures~\ref{fig:num-ellipses}--\ref{fig:num-flower} show the truth and
the reconstructions from noiseless and single-acquisition \(2\%\) data.
The two inclusions remain separated in Case~1. The reconstructions in
Case~2 recover the overall kite outline, with local amplitude variations
near its tips. In Case~3, the outer lobes and central low-value region
are visible in both reconstructions. The viscosity slices show the
change in amplitude in Case~1, rotation in Case~2, and motion of the
high-value region in Case~3.

\clearpage

\subsection*{ Acknowledgments}  
S. Fu is supported by the Fundamental Research Funds for the Central Universities, NPU (Grant No. D5000250416) and by the Natural Science Basic Research Program of Shaanxi (Grant No. 2026JC-YBQN-0003). T. Zhou is partially supported by the National Key Research and Development Program of China (No. 2024YFA1012301), the Zhejiang Provincial Basic Public Welfare Research Program [
Grant Number LDQ24A010001], and NSFC Grant 12371426. 

\subsection*{Data availability statement}
The numerical experiments are based on synthetic data generated using the models and parameters described in this article. The data supporting the numerical results are available from the corresponding author upon reasonable request.
	
\subsection*{Conflict of Interests} 
Hereby we declare there are no conflicts of interest.

\bibliography{references} 
	
\bibliographystyle{alpha}
\end{document}